\documentclass[brochure,english,11pt]{bourbaki}
\usepackage[matrix,arrow]{xy}

\usepackage[pdfpagelabels]{hyperref}

\usepackage{latexsym} 
\usepackage{amsmath}
\usepackage{amssymb}
\usepackage{amsthm}
\usepackage{mathrsfs}
\usepackage{amsfonts}
\usepackage{footnote}%
\usepackage{graphicx}
\usepackage{url}

\usepackage{color}
\usepackage{bbm,verbatim}

\usepackage[latin1]{inputenc}
\usepackage[T1]{fontenc}
\usepackage[frenchb, english]{babel}

\date{Mars 2012}
\bbkannee{64\`eme ann\'ee, 2011-2012}
\bbknumero{1052}

\title{Quantum gravity and the KPZ formula}
\subtitle{after Duplantier-Sheffield}
\author{Christophe GARBAN}
\address{\'Ecole Normale Sup\'erieure de Lyon\\
U.M.P.A.\\
CNRS UMR 5669\\
46 all\'ee d'Italie\\
F--69364 Lyon Cedex 07}
\email{christophe.garban@ens-lyon.fr}

\thanks{Research partially supported by ANR grant MAC2}

\usepackage{multirow} 
\usepackage{array}
\newcolumntype{M}[1]{>{\centering}m{#1}}

\newtheorem{theorem}{Theorem}
\numberwithin{theorem}{section}

\newtheorem{lemma}[theorem]{Lemma}
\newtheorem{proposition}[theorem]{Proposition}

\newtheorem{conjecture}[theorem]{Conjecture}
\newtheorem{question}[theorem]{Question}
\theoremstyle{remark}\newtheorem{definition}[theorem]{Definition}   
\theoremstyle{remark}\newtheorem{remark}[theorem]{Remark}
\theoremstyle{remark}\newtheorem{example}[theorem]{Example}
\def\eqref#1{(\ref{#1})}
\let\qqed=\qed

\def\QED{\qqed\medskip}
\let\qed=\QED

\newcommand{\R}{\mathbb{R}}
\newcommand{\D}{\mathbb{D}}

\newcommand{\C}{\mathbb{C}}
\newcommand{\Z}{\mathbb{Z}}
\newcommand{\N}{\mathbb{N}}

\def\S{\mathbb{S}}

\def \eps {\epsilon}

\def\SLE/{$\mathrm{SLE}$}

\def\Cov#1{\mathrm{Cov}\big[ #1\big]}
\def\Var#1{\mathrm{Var}\big[ #1\big]}
\def \P {{\mathbb P}}
\def \E {{\mathbb E}}

\def\md{\mid}
\def\Bb#1#2{{\def\md{\bigm| }#1\bigl[#2\bigr]}}

\def\Pb{\Bb\P}
\def\Eb{\Bb\E}

\def\Prob#1#2{{\def\md{\bigm| }\P_{#1}\bigl[#2\bigr]}}
\def\<#1{{\langle #1 \rangle}}

\def\M{\mathbf{M}}
\def\m{\mathbf{m}}
\def\K{\mathbb{K}}
\def\re{\overrightarrow{e}}

\def\Zp{\mathcal{Z}}
\def\L{\mathcal{L}}

\def \p {{\partial}}

\def\ising{\mathrm{ising}}

\def\B{\mathcal{B}}
\def\fgam{{\frac {\gamma^2} {2}}}

\def\Sob{\mathcal{H}}

\def\ni{\noindent}
\def\bi{\begin{itemize}}
\def\ei{\end{itemize}}

\def\Ising{MR2680496} 
\def\Onsager{MR0010315}
\def\Kazakov{MR871244}
\def\KPZ{MR947880}

\def\DS{MR2819163}
\def\BS{MR2506765}
\def\GFF{MR2322706}

\def\GCrevisited{MR2642887}
\def\RV{arXiv:0807.1036}
\def\Kahane{MR829798}
\def\GMC{GMC}  

\def\LGu{LGu}
\def\Mu{Mu}
\def\Buz{Buzios}
\def\LGP{MR2438999}
\def\LG{MR2336042}
\def\Msph{MR2399286}
\def\LGM{MR2778796} 
\def\bettinelli{bettinelli}

\def\DupIE{MR1666816}

\def\DupMandel{MR2112128}
\def\DupSaleur{MR889398}
\def\DupKwon{DK88}   

\def\BC{BenjaCurien}  
\def\AC{AngelCurien}  
\def\Omer{MR2024412}

\def\LW{MR1742883}
\def\IEI{MR1879850}
\def\IEII{MR1879851}

\def\lawler{MR2677157}
\def\Dub{MR2525778}
\def\Scott{ArXiv:1012.4797}

\begin{document}

\maketitle

\section{Introduction}

The study of statistical physics models in two dimensions ($d=2$) at their {\it critical point} is in general a significantly hard problem
(not to mention the $d=3$ case). 
In the eighties, three physicists, Knizhnik, Polyakov and Zamolodchikov (KPZ) came up in \cite{\KPZ} with a novel and far-reaching approach in order to understand the critical behavior of these models. Among these, one finds for example random walks, percolation as well as the Ising model. The main underlying idea of their approach 
is to study these models along a two-step procedure as follows:
\bi
\item First of all, instead of considering the model on some regular lattice of the plane (such as $\Z^2$ for example), one defines it instead on a well-chosen ``random planar lattice''. Doing so corresponds to studying the model in its {\it quantum gravity} form. In the case of percolation, the appropriate choice of random lattice matches with the so-called {\bf planar maps} which are currently the subject of an intense activity (see for example \cite{\Buz}). 

\item Then it remains to get back to the actual {\it Euclidean} setup. This is done thanks to the celebrated {\bf KPZ formula} which gives a very precise correspondence between the geometric properties of models in their quantum gravity formulation and their analogs in the Euclidean case.  
\ei

It is fair to say that the nature and the origin of such a powerful correspondence remained rather mysterious for a long time.
In fact, the KPZ formula is still not rigorously established and remains a conjectural correspondence.
The purpose of this survey is to explain how the recent work of Duplantier and Sheffield enables to explain some of the mystery hidden behind this KPZ formula. To summarize their contribution in one sentence, their work implies a beautiful  interpretation of the KPZ correspondence through a uniformization of the random lattice, seen as a Riemann surface.

\medskip
To fully appreciate the results by Duplantier-Sheffield, we will need to introduce beforehand several related concepts and objects. More precisely, the rest of this introduction is divided as follows: first we give a short and informal discussion about {\it quantum gravity}, then we introduce two universality classes of random lattices. Then we will come to the KPZ formula through a specific example (boundary of Random Walks hulls), and finally after stating the main Theorem by Duplantier-Sheffield,  we will state a beautiful conjecture they made.

\subsection{A first glance into quantum gravity}

Quantum gravity is intimately concerned with the following naive question:
\begin{question}\label{NQ}
How does a ``uniformly distributed'' random metric on the sphere $\S^2$ typically look ?
\end{question}
What is naive in this question is the fact that one would first need to specify what we mean by ``a uniform probability measure'' on the space of metrics on $\S^2$. It turns out that defining a natural model of random metric on the sphere $\S^2$ already is a difficult and interesting problem. To illustrate this, let us ask 
a similar naive question in a one-dimensional setting:
\begin{question}\label{QBM}
For any $a,b\in\R^d$, how does a ``uniformly distributed'' path $\gamma: [0,1]\to \R^d$ going from $a$ to $b$ typically look ? 
\end{question}
This naive question was of crucial importance at the time Feynman developed the so-called {\it path integral formulation} of quantum mechanics.

Already in this case, defining properly a ``uniform measure'' on paths was not an easy task. 
Yet, it had been mathematically settled prior to Feynman's work and corresponds to the well-known {\bf Brownian motion}.

In some sense, the purpose of {\it quantum gravity} is to extend Feynman path integrals to  Feynman integrals over {\it surfaces}. Physicists are particularly interested in such an extension, since this would provide a powerful tool to deal with the {\it quantization} of gravitation field theory, a notoriously hard problem. \footnote{This approach towards the quantization of gravitation is called {\it loop quantum gravity}.}
With this background in mind, the problem of defining a proper mathematical object for a ``uniformly chosen random metric on $\S^2$'' thus corresponds to defining a two-dimensional analog of Brownian motion, i.e. a kind of {\it Brownian surface}. 
\medskip

Even though physicists are primarily interested in the above continuum formulation of Question \ref{NQ},
a natural and very fruitful approach is to study an appropriate discretization of it and then to pass to the limit. This brings us to the next subsection.

\subsection{Discretization of Question \ref{NQ} and planar maps}\label{ss.PM}

In the one-dimensional setting, if one asks Question \ref{QBM}, it is not straightforward to come up right away with {\it Brownian motion}. But, if instead we start by discretizing Question \ref{QBM}, say by allowing random $1/n$-steps, then we end up with the model of {\bf random walks}. Brownian motion is then obtained as the scaling limit (as $n\to\infty$) of these rescaled random walks. 

It is thus tempting to apply the same strategy to Question \ref{NQ}, namely to find an appropriate discretization.  Let us explain below a possible discretization which was used initially in the physics literature and was studied extensively recently among the mathematical community. See for example \cite{\Buz} and references therein. We will see in subsection \ref{ss.CM} that there are other ways to discretize Question \ref{NQ} which lead to different universality classes \footnote{Similarly as Random Walks with non-$L^2$ steps converge to other Levy processes than Brownian Motion}, but the discretization below is in some sense the simplest and most natural one regarding the statement of Question \ref{NQ}.

The idea of the discretization we wish to introduce is to consider discrete graphs, with say $n$ faces, which have the topology of a sphere $\S^2$ and for which the metric $\rho_n$ will correspond (up to a rescaling factor) to the graph distance, i.e. $\rho_n:=n^{-a}d_{gr}$ for some exponent $a>0$. The exponent $a$ will need to be well chosen  as in the case of Random Walks, where space needs to be rescaled by $\sqrt{n}$ in order to obtain a limit. If we define our discretization in such a way that for each $n\geq 1$, there are finitely many  such graphs, we can pick one uniformly at random (in the spirit of Question \ref{NQ}) which thus gives us a random metric space $(M_n,\rho_n)$. We can then ask the question of the scaling limit of these random variables $(M_n,\rho_n)$ as $n\to \infty$ in the space $(\K,d_{GH})$ of all isometry classes of compact metric spaces, endowed with the Gromov-Hausdorff distance $d_{GH}$
\footnote{Informally, if $(E_1,d_1)$ and $(E_2,d_2)$ are two compact metric spaces, then $d_{GH}(E_1,E_2)$ is computed as follows: we embed $E_1$ and $E_2$ isometrically into some larger metric space $(E,d)$ and we compute using the common distance $d$ the distance in the Hausdorff sense between the two embeddings. Then we take the infimum over the possible such embeddings. See \cite{\Buz}.}.
The advantage of this setup is that $(\K,d_{GH})$ is a complete, separable, metric space (a Polish space) and is thus suitable to the analysis of convergences in law and so on.  
Note here, that even if one could prove that $(M_n,\rho_n)$ converges to a limiting random object $(M_\infty,\rho_\infty)$, it is not clear a priori that the topology is preserved at the scaling limit or in other words, it needs to be proved whether $(M_\infty,\rho_\infty)$ a.s. has the topology of a sphere or not.  If all these steps can be carried on, then this would give us a good candidate $(M_\infty,\rho_\infty)$ for the random object used in Question \ref{NQ}. 

\begin{figure}[!htp]
\begin{center}
\includegraphics[width=0.65\textwidth]{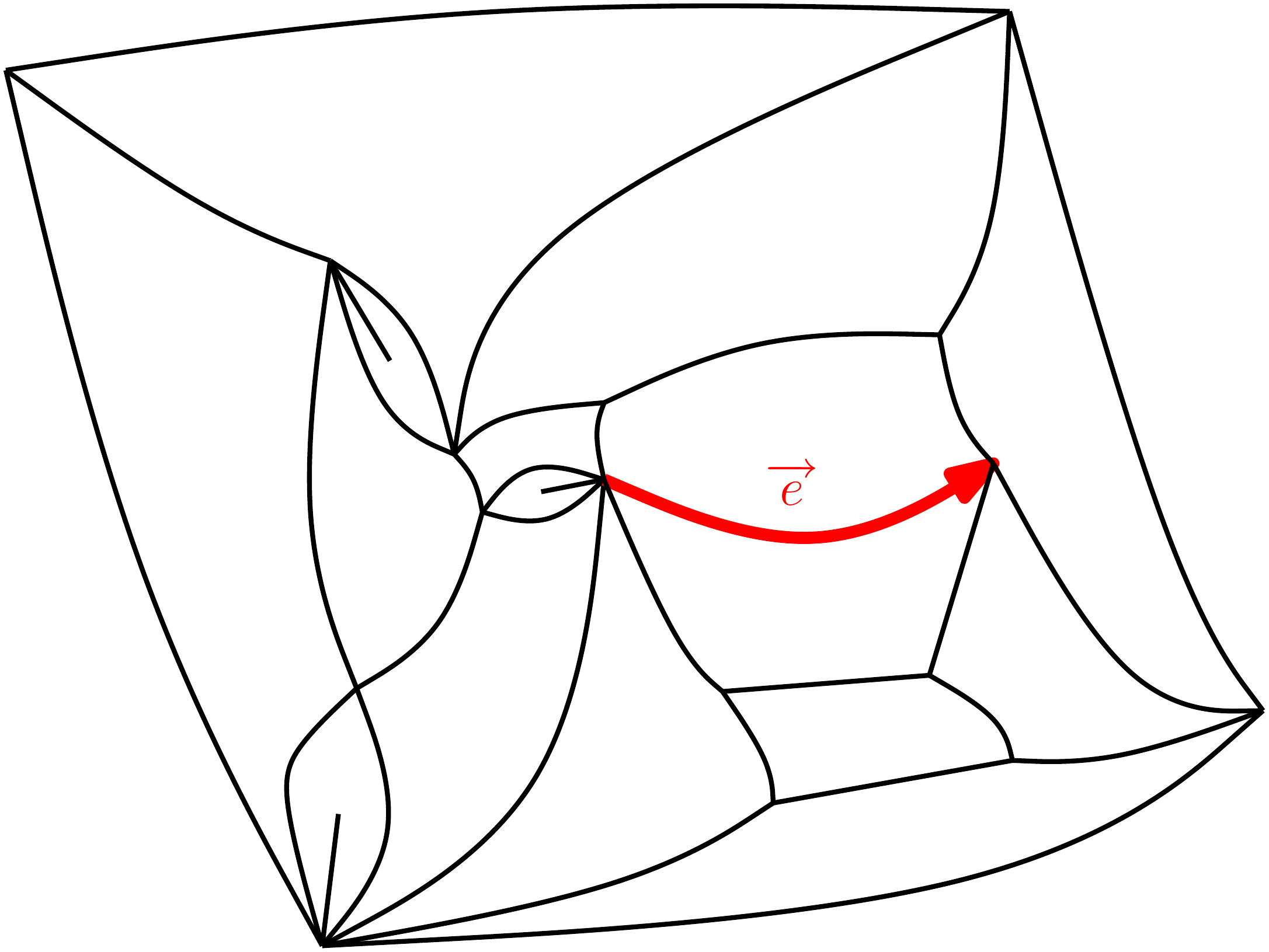}
\end{center}
\caption{This is a planar map of the sphere $\S^2$ with exactly 17 squares (this includes the exterior square which is also in the sphere).}\label{f.quadrangulation}
\end{figure}

\medskip
Let us now introduce one specific discretization.
\begin{definition}[planar map, following \cite{\LGu}]\label{d.pm}
A {\bf planar map} $M$ is a proper embedding of a finite and connected graph into the two-dimensional sphere $\S^2$, which is viewed up to orientation preserving homeomorphisms of $\S^2$ (i.e. up to ``deformations''). 
Loops and multiple edges are allowed. The faces of $M$ are identified with the connected components of $\S^2\setminus M$ and the degree of a face $\mathbf{f}$ is defined as the number of edges incident to $\mathbf{f}$, with the additional rule that if both sides of an edge belong to the same face, this edge is counted twice. 

Finally, for combinatorial reasons, it is often convenient to consider {\bf rooted} planar maps, meaning that one particular oriented edge $\re$ is distinguished. The origin  of that root edge $\re$ is called the {\bf root vertex} $\emptyset$. See figure \ref{f.quadrangulation} for an instance of a planar map where all faces happen to be squares.  
\end{definition}

\begin{definition}[$p$-angulations of the sphere]\label{d.pa}
For any integer $p\geq 3$, let $\M_n^p$ be the set of all rooted planar maps with $n$ faces, where each face has degree $p$. The elements of $\M_n^p$ are called {\bf rooted planar $p$-angulations}. ($p=3$ corresponds to triangulations and $p=4$ to quadrangulations).
\end{definition}

Since planar maps are defined up to deformations, there are finitely many rooted planar maps in $\M_n^p$, for any $n\geq 1$. For example, one has 
\begin{equation*}\label{}
\# \M_n^4 = \frac{2}{n+2} 3^n \mathrm{Cat}_n = \frac{2\cdot 3^n}{(n+2)(n+1)} \binom{2n}{n}\,.
\end{equation*}
The appearance of the $n$-th Catalan number here is explained through the celebrated  {\bf Cori-Vauquelin-Schaeffer bijection}, which gives a one-to-one correspondence between labelled plane trees and rooted quadrangulations of the sphere (see for example \cite{\Buz}). 

\medskip 

For a fixed integer $p\geq 3$, we will denote by $\m_n\in \M_n^p$ a sample of a planar map uniformly distributed over $\M_n^p$. 
Since in Question \ref{NQ}, we were interested in a uniformly distributed random metric on $\S^2$, it makes sense to consider the random variable $(\m_n, n^{-a}d_{gr})$ seen as a random point in the above space $(\K,d_{GH})$.

The study of these random planar maps has now a long history (see for example \cite{\Buz, \LGP,\LG,\Msph}) and has culminated in the following breakthrough result proved independently by Le Gall (\cite{\LGu}) and Miermont (\cite{\Mu})
\footnote{To be more precise, Miermont's proof is restricted to the case $p=4$ but gives slightly stronger estimates on the structure of the geodesics while Le Gall's proof is more general in the sense that it proves the result also for $p=3$ and all even $p$ greater than 4.}:
\begin{theorem}[Uniqueness of the scaling limit]\label{th.u}
There exists a (unique) random compact metric space $(\m_\infty,D^*)$ with values in $\K$ such that for any $p\in 3\cup (2\N +4)$, there is a positive constant $\lambda_p$ such that as $n\to \infty$, 
\begin{equation*}\label{}
(\m_n, \frac {\lambda_p} {n^{1/4}}\, d_{gr}) \overset{(d)}{\longrightarrow} (\m_\infty, D^*)\,,
\end{equation*}
in the Gromov-Hausdorff sense. This random compact metric space is called the {\bf Brownian map}. 
\end{theorem}

Furthermore the following property holds. (It was first established in \cite{\LGP}. See also \cite{\Msph} for a different proof.)
\begin{theorem}[Sphericity of the Brownian map]\label{th.s}
Almost surely, the random metric space $(\m_\infty, D^*)$ is homeomorphic to the sphere $\S^2$.
\end{theorem}

This last theorem thus ensures that in some sense the Brownian map $(\m_\infty, D^*)$ is a good candidate for the ``uniform'' random metric on $\S^2$ considered in Question \ref{NQ}. The only drawback being that the metric space $(\m_\infty, D^*)$ which is a.s. homeomorphic to $(\S^2, \|\cdot\|_{\R^3})$ is not provided with a ``canonical'' embedding into the sphere
\footnote{This would be the case for example if $(m_\infty, D^*)$ happened to be a nice Riemann surface, but it is not.}. We will come back to this question of  embedding later.

\subsection{Quantum gravity coupled with ``matter''}\label{ss.CM}
Coming back to our earlier motivation, i.e. the study of statistical physics models on regular lattices, 
we will now introduce models of random lattices which are naturally associated with statistical physics models. The first one will be associated to the Ising model and will fall into a different {\it quantum gravity} universality class than the above planar maps.

\subsubsection{Quantum gravity coupled with Ising model}

The standard and simplest way in statistical physics to model {\it ferromagnetic} matter (like a piece of iron for instance) is through the so-called {\bf Ising model}. It is defined on a graph $G=(V,E)$ which is supposed to represent the metallic structure of our ferromagnet. For a flat piece of iron, one might choose $\Z^2$ for example. On a finite graph $G=(V,E)$, it is defined as a probability measure on {\bf spin configurations} $\sigma\in \{-1,1\}^V$, where if $x\in V$, $\sigma_x=+1$ means that the {\it spin of the atom} at site $x$ is oriented in the upward direction. The probability measure $\P=\P_J$ is such that each spin configuration $\sigma$ has a probability proportional to 
\begin{equation*}\label{}
\Prob{J}{\sigma}\propto \exp  \Big( \;J\sum_{e= \<{x,y} } \sigma_x \sigma_y \Big)\,,
\end{equation*}
where the parameter $J$ represents the strength of the electromagnetic interaction between atoms
\footnote{Often, instead of the parameter $J$, one uses the {\bf inverse temperature} $\beta:=\frac 1 {k_B T}$ 
and the Ising model is defined so that $\Prob{\beta}{\sigma}\propto e^{\;\beta\sum_{e= \<{x,y} } \sigma_x \sigma_y} $, but mainly for notational reasons, we will not use this point of view here.}
.
Note that the higher $J$ is, the greater the tendency is for the spins to align in the same direction. 
In order to express the probability measure $\P$ more explicitly, it is natural to introduce the {\bf partition function} of the graph $G$
\begin{equation*}\label{}
Z_J=Z_J(G):=\sum_{\sigma\in \{-1,1\}^V} \exp \Big(\, J \sum_{e=\<{x,y}} \sigma_x\sigma_y \Big)\,,
\end{equation*}
so that 
\begin{equation*}\label{}
\Prob{J}{\sigma}=\frac 1 {Z_J} \, \exp \Big( \;J\sum_{e= \<{x,y} } \sigma_x \sigma_y \Big)\,.
\end{equation*}

The Ising model has been studied extensively in the physics and mathematical community over the last 70 years. See for example \cite{\Ising} and references therein for a recent breakthrough paper on the subject.

Now the idea of studying such a model from a quantum gravity perspective is to sample a spin configuration $\sigma$ together with a random base-graph $G$. This is what one calls {\it coupling}  quantum gravity with the Ising model (i.e. we are looking for a probability measure on pairs $\{ G=(V,E),\sigma\in\{-1,1\}^V \}$ ). One way to proceed would be to first sample the random lattice $G$, using the above planar maps $\m_n$ uniformly chosen among $\M_n^p$ and then to sample an Ising configuration $\sigma$ on the graph  $\m_n$ according to the above model $\P=\P_{J,\m_n}$. Doing so would not correspond so much to {\it coupling} quantum gravity with the Ising model, since the marginal of $(\m_n,\sigma)$ on the graph-component would just be the above model of planar graphs. So we will use a more intricate way of coupling the graph with its spin-configuration. 

\begin{figure}[!htp]
\begin{center}
\includegraphics[width=0.6\textwidth]{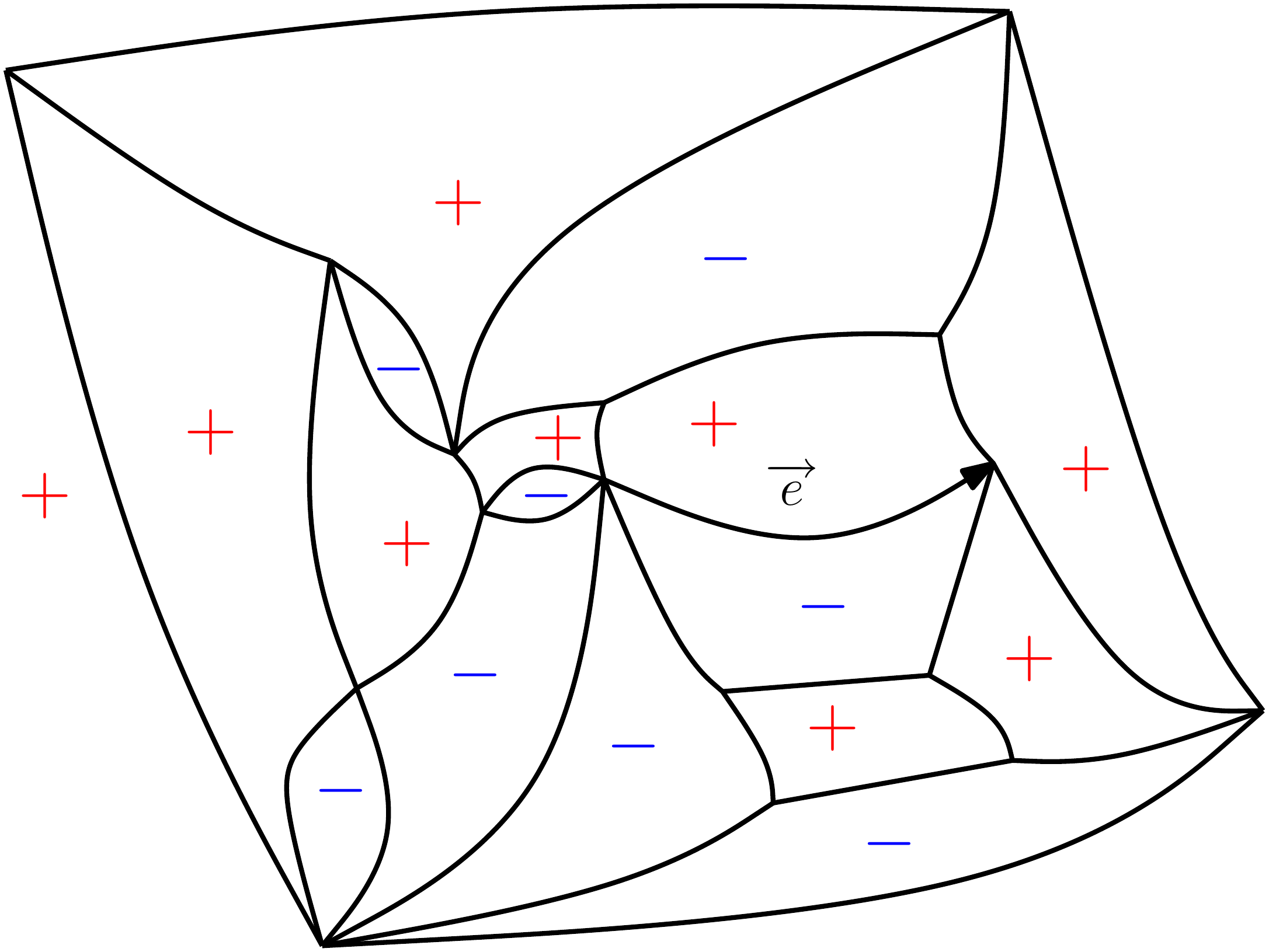}
\end{center}
\caption{This is an Ising spin configuration on the same quadrangulation as in Figure \ref{f.quadrangulation} (were we removed  the inside edges for convenience).}\label{f.ising}
\end{figure}

If we restrict ourselves to the case of planar graphs, the standard quantum gravity/Ising coupling works as follows: fix $p\geq 3$ and some large $n\geq 1$.  For integrability reasons, the spins will be indexed by the faces of the planar maps $\m\in \M_n^p$ instead of the vertices. (In other words, we will consider our Ising model on the dual graphs $\m^*$ which are no longer in $\M_n^p$). With this slight change, a natural idea is to define a measure $\P=\P_{J,n}$ on couplings $(\m_n,\sigma)\in \M_n^p \times \{-1,1\}^n$ (recall that planar maps in $\M_n^p$ have exactly  $n$ faces) so that 
\begin{equation*}\label{}
\Prob{J,n}{\m,\sigma}\propto \exp \Big(  J \sum_{
f\sim f' 
} \sigma_f \sigma_{f'} \Big)\,,
\end{equation*}
where the sum is over pairs of adjacent faces in $\m$. See Figure \ref{f.ising}. 

Note that, unless $J\equiv 0$, the marginal on $\m$ (the graph-component) is no longer uniformly chosen among $\M_n^p$, but follows a law which highly depends on the magnetic strength $J$. In probabilistic terms, the marginal on $\m_n$ corresponds to weighting $\m_n$ by its Ising partition function $Z_J(\m_n)$. As in the case of a fixed graph $G$, it is natural to define the following {\it annealed} partition function 
\begin{equation*}\label{}
\Zp_n(J):=\sum_{\m\in \M_n^p, \sigma\in \{-1,1\}^\m} \exp \big( \, J \sum \sigma_f \sigma_{f'} \big)=\sum_{\m\in \M_n^p} Z_J(\m_n)
\end{equation*}

Now, recall the KPZ formula is supposed to relate a {\bf critical} Euclidean Ising model (say on $\Z^2$) with a {\bf critical} Ising model 
defined together with a random planar map $\m\in \bigcup_n \M_n^p$ in the above fashion. Therefore, we need a way to detect a phase transition in each model, Euclidean and quantum gravity. Among the many ways to ``detect'' phase transitions, one is to notice a  failure of analyticity for certain thermodynamic quantities. In our case, partition functions enable us to detect a phase transition. More precisely, for the Euclidean Ising model on $\Z^2$, for each $n\geq 1$, let $Z_n(J)$ be the above partition function for the graph $G_n=[-n,n]^2\cap \Z^2$. Then, from the work of Onsager \cite{\Onsager}, it can be shown that as $n$ goes to infinity, the functions $J\mapsto \frac 1 {n^2} \log Z_n(J)$ converge to a limiting function $f(J)$, called the {\it free energy} which is analytic for all values of $J > 0$ except at $J=J_c^{\mathrm{eucl}}= \frac 1 2 \log(1+\sqrt{2})$. 
On the quantum gravity side, the natural quantity to look at is the so called {\bf grand canonical ensemble} partition function defined as 
\begin{eqnarray*}\label{}
\Zp(\beta,J) & := & \sum_{n\geq 1}\sum_{\m\in \M_n^p, \sigma \in \{-1,1\}^n} e^{-\beta |\m|} e^{J\sum_{f\sim f'} \sigma_f \sigma_{f'}} \\
&=& \sum_{n\geq 1} e^{-n \beta} \Zp_n(J)\,.
\end{eqnarray*}
The reason why such an annealed partition function is considered is because it was computed exactly by Kazakov (\cite{\Kazakov})
in the case of quadrangulations ($p=4$). His computation relied on a deep relationship with some (random) matrix models.
This type of exact computation occurs for what one calls {\it exactly solvable models} and is the reason why we changed slightly the definition of Ising model from vertices to faces. The exact computation by Kazakov enables to detect a bi-critical point $(\beta_c,J_c)$
\footnote{When $p=4$, one has $J_c=J_c^{QG}=\ln 2$} around which the analyticity of $\Zp(\beta,J)$ is broken.  
The analysis of $\Zp$ around its bi-critical point gives detailed information on the critical behavior of the system in its quantum gravity formulation. See \cite{\LGM} and Appendix B in \cite{\DupMandel} for thorough discussions of this.

As mentioned above, the scaling limit as $n\to\infty$ of these planar maps $\m_n\in \M_n^4$ weighted by $Z_{J_c=\ln 2}(\m_n)$ will fall in a different universality class than the ``uniform'' (or unweighted) planar maps.  We will come back to this later.

\subsubsection{Quantum gravity coupled with random walks, self avoiding walks or percolation}

If one wants to study statistical models such as random walk, percolation or self-avoiding walks (SAW), it turns out that all of them are naturally coupled with the planar maps introduced in subsection \ref{ss.PM}. This regime corresponds to what one calls the {\bf pure gravity} regime. It means in some sense that the random geometry of the planar map is insensitive to the model it is coupled to (as opposed to what happens with Ising model). We will describe this pure gravity case in more details in the next subsection through the example of random walk.

\subsection{The KPZ formula in the special case of random walk and Brownian motion}

\begin{figure}[!htp]
\begin{center}
\includegraphics[width=\textwidth]{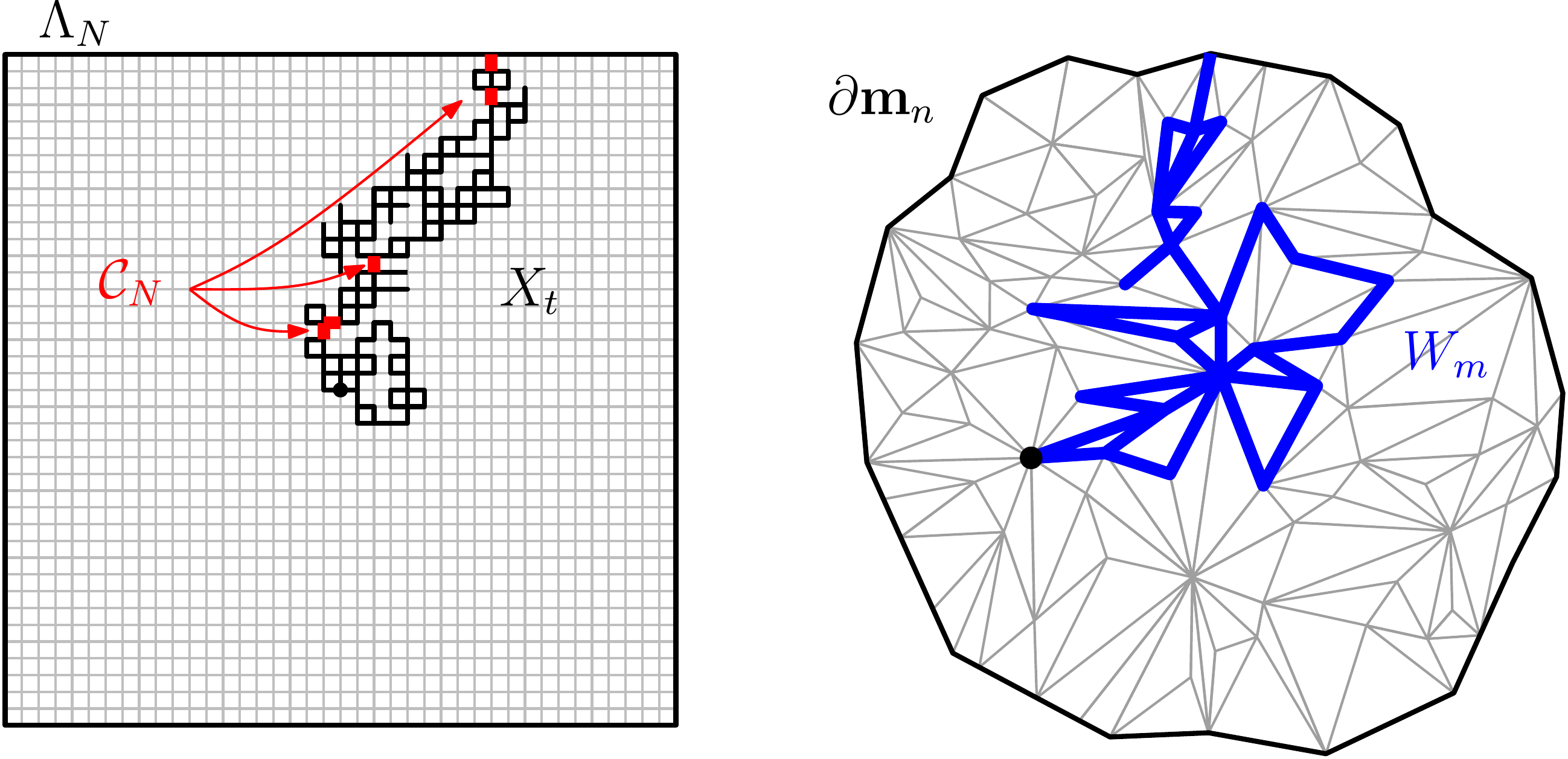}
\end{center}
\caption{}\label{f.euclidVSquantum}
\end{figure}

Let us start by the Euclidean case. Consider a random walk $X_t$ on the domain $\Lambda_N:=[-N,N]^2$ starting at the origin until it reaches the boundary $\p \Lambda_N$. One can look at several subsets of interest about this random walk among which:
\bi
\item[(i)] The {\it range} $\mathcal{R}_N$ of $X_t$, i.e. $\mathcal{R}_N:=\{X_t\}$.
\item[(ii)] The set of {\it cut-points} $\mathcal{C}_N$ of $X_t$, i.e. the set of points $x\in \Lambda_N$ so that removing $x$ disconnects the range into 2 disjoint components.
\item[(iii)] The set of {\it frontier points} $\mathcal{F}_N$, i.e. all the points on the range that are connected to $\p \Lambda_N$ in the complement $\Lambda_N \setminus \mathcal{R}_N$. 
\ei
See figure \ref{f.euclidVSquantum}. 
One way to measure the typical size of a random subset of $\Lambda_N$ is via the following notion of scaling exponent. 
If $K=(K_N)_{N\geq 1}$ is a certain sequence of random subsets of $\Lambda_N$, then its {\bf Euclidean scaling exponent}
is defined as
\begin{equation}\label{e.ese}
x=x(K):=\lim_{N\to \infty} \frac{\log \Eb{|K_N|/N^2}}{\log 1/N^2}\,.
\end{equation}
Note that we are being informal here since it might be that this limit does not exist. To be more rigorous, one should define instead upper and (lower) scaling exponents by using $\limsup$ (and $\liminf$) instead. Nevertheless, to simplify the exposition we will neglect this issue in the remaining of this text and will assume that the random sets we will consider are such that these limits always exist (proving such a convergence can be very hard for some models and remains in many cases an open problem, for example scaling exponents for critical percolation on $\Z^2$). 
Note also that the greater the exponent $x=x(K)\in[0,1]$ is, the ``smaller'' the random sets $K_N$ are since having scaling exponent $x$ means that asymptotically, $K_N$ is of size $\approx N^{2-2x}$.

The {\bf scaling exponents} for the above subsets of random walk paths have been the focus of an intense activity over the last 20 years or so. In fact the road which culminated in their evaluation is quite an interesting story. 
Their values were first conjectured by Duplantier and Kwon in \cite{\DupKwon} based on conformal invariance and numerical simulations. Somewhat surprisingly, the powerful machinery  of {\it conformal field theory} was unable at giving theoretical predictions for these exponents. 
In 1998, Bertrand Duplantier obtained in \cite{\DupIE} (on a non-rigorous basis) these exponents with a completely different approach: his idea was to study a random walk path on a random planar map $\m_n$ and to compute the scaling exponents, called {\bf quantum scaling exponents} for the frontier and cut points of the random walk on $\m_n$. Once this was achieved, he was able to recover his previous conjecture via the KPZ formula.
His work was partly based on the seminal work \cite{\LW} on Brownian intersection exponents. Finally, the story ends with the 
mathematical derivation of these exponents in \cite{\IEI,\IEII} which used yet another approach: the Schramm-Loewner-Evolution  $SLE_\kappa$ processes (with $\kappa=6$ here).

Here is more precisely how it goes. Imagine that we have at our disposal a natural model of planar maps $\m_n$ with boundary (i.e. with the topology of a disk), instead of planar maps $\m_n$ with the topology of a sphere. This can be done in such a way that it gives a different, yet interesting model of planar maps:  see \cite{\bettinelli} where Bettinelli considers random quadrangulations $\m_n$ with $n$ faces and which are such that $\p \m_n$ is of order $\sqrt{n}$. Let thus $\m_n$ be sampled uniformly among rooted quadrangulations with $n$ faces and boundary of length $\sqrt{n}$.
Let $W_m$ be a {\bf simple random walk} on this graph $\m_n$ starting at the root of $\m_n$ until it reaches the boundary $\p \m_n$. See figure \ref{f.euclidVSquantum}.
One can consider the exact analogs of the above sets, i.e:
\bi
\item[(i)] The {\it range} $\mathcal{R}_n$ of $W_m$, i.e. $\mathcal{R}_n:=\{W_m\}$.
\item[(ii)] The {\it cut-points} $\mathcal{C}_n$ of $W_m$, i.e. the set of points $x\in \m_n$ so that removing $x$ disconnects the range into 2 disjoint components.
\item[(iii)] The frontier points $\mathcal{F}_n$, i.e. all the points on the range that can be connected to $\p \m_n$ via the complement $\m_n \setminus \mathcal{R}_n$. 
\ei

Similarly to the Euclidean case, if $K=(K_n)_{n\geq 1}$ is a sequence of random subsets of 
$\m_n$, then its {\bf quantum scaling exponent} is defined as:

\begin{equation}\label{e.qse}
\Delta=\Delta(K):=\lim_{n\to \infty} \frac{\Eb{ \log |K_n| / |\m_n| }}{\log 1/ |\m_n|} = \lim_{n\to \infty} \frac{\Eb {\log |K_n| /n} }{\log 1/n}\,,
\end{equation}
where, once again, we implicitly assume here that the limits exist. 

Even though it was not exactly in the same setup, Duplantier determined in \cite{\DupIE} the values of these quantum scaling exponents. 
He found  $\Delta(\mathcal{F})=1/2$ and $\Delta(\mathcal{C})=3/4$. Now for this universality class (i.e. the {\bf pure gravity} case of ``uniform'' planar maps $\m_n$), the KPZ formula reads as follows: the scaling exponents $x$ and $\Delta$ of a certain subset of a random walk path in the Euclidean and quantum cases are related by the following quadratic expression:
\begin{equation}\label{KPZ1}
x=\frac{2}3 \Delta^2 + \frac 1 3 \Delta\,.
\end{equation}
One thus finds $x(\mathcal{F})=1/3$ and $x(\mathcal{C})=5/8$. Note that for the {\it range} $\mathcal{R}$, it is well known that an Euclidean  random walk $X_t$ visits about $N^2/\log N$ sites before touching $\p \Lambda_N$, this means that $x(\mathcal{R})=0$ which is a fixed point of the above KPZ formula. This thus suggests that $\Delta(\mathcal{R})=0$ as well, which in turn means that  the random walk $W_m$ should visit $|V(\m_n)|^{1-o(1)}=n^{1-o(1)}$ sites of the planar map $\m_n$. This is a non-trivial fact in its own, which illustrates that the KPZ formula (if true) can be powerful in both directions.

\subsection{The KPZ formula in general}

The planar maps $\m_n$ uniformly chosen among some $\M_n^p$ for some $p\geq 3$ are the natural universality class for studying random walk on the quantum gravity side. Nevertheless, one might also consider an independent random walk on a planar map $\m_n$ sampled together with an Ising model as we explained above. If one would do so, this would affect the typical size of the sets $\mathcal{C}_n$ and $\mathcal{F}_n$ and one would find different values for the quantum scaling exponents $\Delta_\ising(\mathcal{F})$ and $\Delta_\ising(\mathcal{C})$. This reflects the fact that the two different procedures we gave to sample planar maps do not fall asymptotically in the same universality class. For this Ising universality class, the KPZ relation reads as follows:
\begin{equation}\label{KPZising}
x=\frac{3}4 \Delta_\ising^2 + \frac 1 4 \Delta_\ising\,.
\end{equation}

Such a correspondence can be very useful to relate critical exponents for the Ising model in its quantum gravity form with the analogous Euclidean critical exponents. For example, one might hope that if $K_n$ denotes the largest $+$ cluster in the random map $\m_n$ weighted according to $Z_J(\m_n)$, then its quantum scaling exponent $\Delta_\ising$ should relate to its Euclidean analog through the above formula. Yet, one has to be careful with such statements since, especially away from the {\it pure gravity case} (i.e. uniform planar maps), it is in general a subtle affair to know which critical exponents are allowed to go through the KPZ formula or not. 
\footnote{The subtlety lies partly in the fact that one needs the set $K$ to be independent of the field $e^{\gamma h}$ in the main Theorem \ref{th.main}, while here the largest $+$ cluster would be highly correlated with the field $e^{\gamma h}$.} 
In the case of the Ising model, quantum critical exponents are usually derived from the analytic study of $\Zp(\beta,J)$ near its bi-critical point and are then transfered to the Euclidean setting using~\eqref{KPZising}. See \cite{\DupMandel}.

One can now state the general KPZ formula: consider some statistical physics model defined together with a random planar lattice (in the suitable universality class). Then, the quantum and Euclidean scaling exponents $\Delta$ and $x$ which describe the critical properties of the model, are such that  the following quadratic KPZ formula holds
\begin{equation}\label{KPZ}
x=\frac{\gamma^2} 4 \Delta^2 + (1-\frac {\gamma^2} 4) \Delta\,,
\end{equation}
where the parameter $\gamma\in\R_+$ determines in which universality class the model is. We thus find that {\it pure gravity} 
corresponds to $\gamma=\sqrt{8/3}$ while Ising model corresponds to  $\gamma=\sqrt{3}$. More generally, it follows from this correspondence that critical two-dimensional models form in some sense a one-dimensional family space. This is consistent with the recent theory of $SLE_\kappa$ processes which are aimed at describing interfaces of critical two-dimensional systems and also form a one-dimensional family of processes indexed by $\kappa\geq 0$. In fact, the two parameters are related one to another by $\gamma\equiv \sqrt{\kappa}$.

\subsection{Why go through quantum gravity ?}

The KPZ formula becomes particularly useful (at least on a non-rigorous level) when the critical exponents of a particular model happen to be much easier to compute, say, in its quantum gravity form than in its Euclidean one (or vice-versa). This is exactly what happens in the case of random walk (\cite{\DupIE}). 

\subsubsection{Quantum gravity coupled with random walk}.
Imagine we want to study the exterior frontier $\mathcal{F}_t$ of a random walk $X_t$ in $\Z^2$ as $t$ increases ($\mathcal{F}_t$ is defined here as the set of points in $\mathcal{R}_t$ which are connected to infinity in $\Z^2\setminus \mathcal{R}_t$). For each $T\geq 1$, the evolution of $\mathcal{F}_t$ for $t\geq T$ will not depend of what the random walk did inside the hull $\mathcal{D}_T$, defined as the complement of the unique infinite connected component of $\Z^2 \setminus \mathcal{R}_T$. Yet the evolution of $\mathcal{F}_t, t\geq T$ will depend of the complicated law of the boundary $\p \mathcal{D}_T=\mathcal{F}_T\subset \Z^2$ which as $T\to \infty$ looks more and more like a fractal set. On the quantum gravity side, such geometric processes as $\mathcal{F}_t$ behave in some sense in a much more Markovian way.

To illustrate this, let us briefly describe a striking recent work by Benjamini and Curien \cite{\BC}. In their work, they rely on an infinite version $Q_\infty$ of the above rooted quadrangulations of the sphere $\m_n\in \M_n^4$. $Q_\infty$ is called a rooted UIPQ (for Uniform Infinite Planar Quadrangulations) and can be seen as a local limit of the quadrangulations $\m_n$ as $n\to \infty$.  They consider a random walk $W_m$ starting at the root of $Q_\infty$ and they explore the planar map $Q_\infty$ along the random walk path $W_m$. The crucial observation which goes back to \cite{\Omer} is that if at time $m$, $\mathcal{F}_m$ denotes the exterior boundary of the domain which was explored by the random walk so far \footnote{some care is needed here, since one has to prove that under the law of $Q_\infty$, there is a unique infinite connected component in $Q_\infty \setminus \mathcal{R}_m$, but we will not enter in more details here. See \cite{\BC}.}, then the law of the infinite connected component $Q_\infty \setminus \mathcal{D}_m$ depends only on the length of $\p \mathcal{D}_m = \mathcal{F}_m$ and not at all on its shape.
This simplifies things tremendously. Yet, it would remain to control how $|\mathcal{F}_m|$ and $|\mathcal{D}_m|$ both increase as a function of $m$. This part turns out to be difficult but Benjamini and Curien are able to obtain in \cite{\BC} that as $m\to \infty$, $|\mathcal{F}_m|$ behaves like 
$|\mathcal{D}_m|^{1/2}$. This provides a rigorous derivation of the identity  $\Delta(\mathcal{F})=1/2$. (They also obtain another quantum scaling exponent for the so-called {\it frontier points}). 

These nice spatial Markovian properties were also crucial in the seminal work on the quantum scaling exponents of random walk by Duplantier (\cite{\DupIE}), where various quantum scaling exponents for the random walk were determined through an asymptotic analysis of the partition function of random walks coupled with planar maps. 


\subsubsection{Quantum gravity coupled with percolation}
Critical percolation has also been successfully analyzed from the quantum gravity perspective. Various quantum critical exponents for percolation were determined in \cite{\DupSaleur} and were translated to the Euclidean setting via the KPZ formula. This is probably the first explicit use in the physics literature of the KPZ formula. 
   In the mathematics literature, there is a work in progress by Angel and Curien \cite{\AC}, which gives among other things a rigorous determination of the critical points for bond percolation on the UIPT and UIPQ (the first one is the analog of the above $Q_\infty$ with triangles instead). Using the Markovian structure we briefly sketched above, they are able to compute some of the quantum scaling exponents of these critical percolation models, which is a very interesting step towards the understanding of statistical physics models on random lattices.

\subsubsection{Quantum gravity coupled with Ising model}
The reason why it is interesting to study the Ising model in its quantum gravity form lies in the fact that (as was mentioned above), Kazakov was able to compute exactly the annealed partition function 
\begin{equation*}\label{}
\Zp(\beta,J)=\sum_{n\geq 1} e^{-\beta \, n} \sum_{\m\in \M_n^4} Z_J(\m_n)\,.
\end{equation*}
As hinted previously, various quantum critical exponents for the Ising model can be extracted from the behavior of $\Zp(\beta,J)$ near $(\beta_c,J_c)$ and can be translated to the Euclidean setting via the KPZ formula~\eqref{KPZ} with $\gamma=\sqrt{3}$. 

It is now time to introduce the main result we wish to explain in this survey.

\subsection{The main result by Duplantier-Sheffield and a conjecture on the embedding of planar maps}

What remains mysterious so far in the KPZ correspondence we just described is the fact that it relates scaling exponents of sets 
which do not ``live'' on the same space. For example it relates the scaling exponent of cut-points of a random walk living in $\Lambda_N$ with the (quantum)-scaling exponent of cut-points of a random walk living on a random planar map $\m_n$. 
To overcome this, Duplantier and Sheffield discovered a setup in which the random set we are interested in, say $K$, lives on some single space (say $[0,1]^2$ or $\S^2$), but its size, given in terms of its scaling exponent, can be measured in two different ways: an ``Euclidean'' way which will output some Euclidean scaling exponent $x=x(K)$ and a ``quantum'' way which will output a possibly different quantum scaling exponent $\Delta=\Delta(K)$. Their setup is built so that $x(K)$ and $\Delta(K)$ will satisfy to the general KPZ relation~\eqref{KPZ}. (More precisely for each value of $\gamma\geq0$, they found an adequate setup). 

Before going into the details of their setup, one sees here that we are getting closer to our initial naive Question \ref{NQ}. Indeed, if one had at our disposal a ``uniform'' random metric $\rho$ on $\S^2$, one could evaluate the size of $K$ either using the Euclidean metric $\|\cdot \|$ (which would give us a scaling exponent $x$) or using the random metric $\rho$ (which would give us an exponent $\Delta$). Now from the preceding discussions, $x$ and $\Delta$ should satisfy to~\eqref{KPZ} with $\gamma=\sqrt{8/3}$. Unfortunately we did not quite succeed yet in answering Question \ref{NQ}. Indeed, Theorems \ref{th.u} and \ref{th.s} provide us with a ``uniform'' probability distribution on compact metric spaces which a.s. have the topology of a sphere, but 
there is (at least for the moment) no canonical way to embed the limiting space $(\m_\infty, D^*)$ into the sphere. In fact, because of this embedding issue, a proper answer to Question \ref{NQ} remains an important open problem. 

Yet, before passing to the limit $n\to \infty$, there are several ways to ``naturally'' embed the planar maps $\m_n$ into the sphere. We present two of them.

\subsubsection{Embedding of planar maps seen as Riemann Surfaces}

The main idea here is that one can view each planar map $\m_n\in \M_n^4$ as a Riemann surface. For this, view each face as a polygon (here a square) on which we give the obvious flat conformal structure given formally by $z\mapsto z$. 
By the Schwarz reflexion principle, one can easily glue together the conformal structures of two adjacent faces along their edge. Some care is needed around each vertex $x\in \m_n$ since the angle around $x$ might not be $2\pi$. Yet, one can use a local chart around each such $x$ of the form $z \mapsto z^{\frac 4 k }$ where $k$ is the degree of vertex $x$ (this corresponds to the fact that conic singularities of complex manifolds are {\it removable}). 
Altogether, this gives us a complex manifold of dimension one endowed with a finite atlas indexed by the set of edges and vertices of $\m_n$. 
By the Riemann uniformization theorem, since $\m_n\simeq \S^2$, $\m_n$ equipped with the above conformal structure can be mapped conformally to $\S^2$ and the embedding is {\bf unique} up to M\"obius maps of the sphere $\S^2$.  We thus found a natural way (up to M\"obius transformations) to embed any planar maps $\m_n$ into $\S^2$. The same idea enables us to embed in a conformal way planar maps $\m_n$ with the topology of a disc (with $\p \m_n\neq \emptyset$) into the disc $\D$ or into $[0,1]^2$. 

\subsubsection{Embedding of planar maps via circle packings}
\begin{figure}[!htp]
\begin{center}
\includegraphics[width=\textwidth]{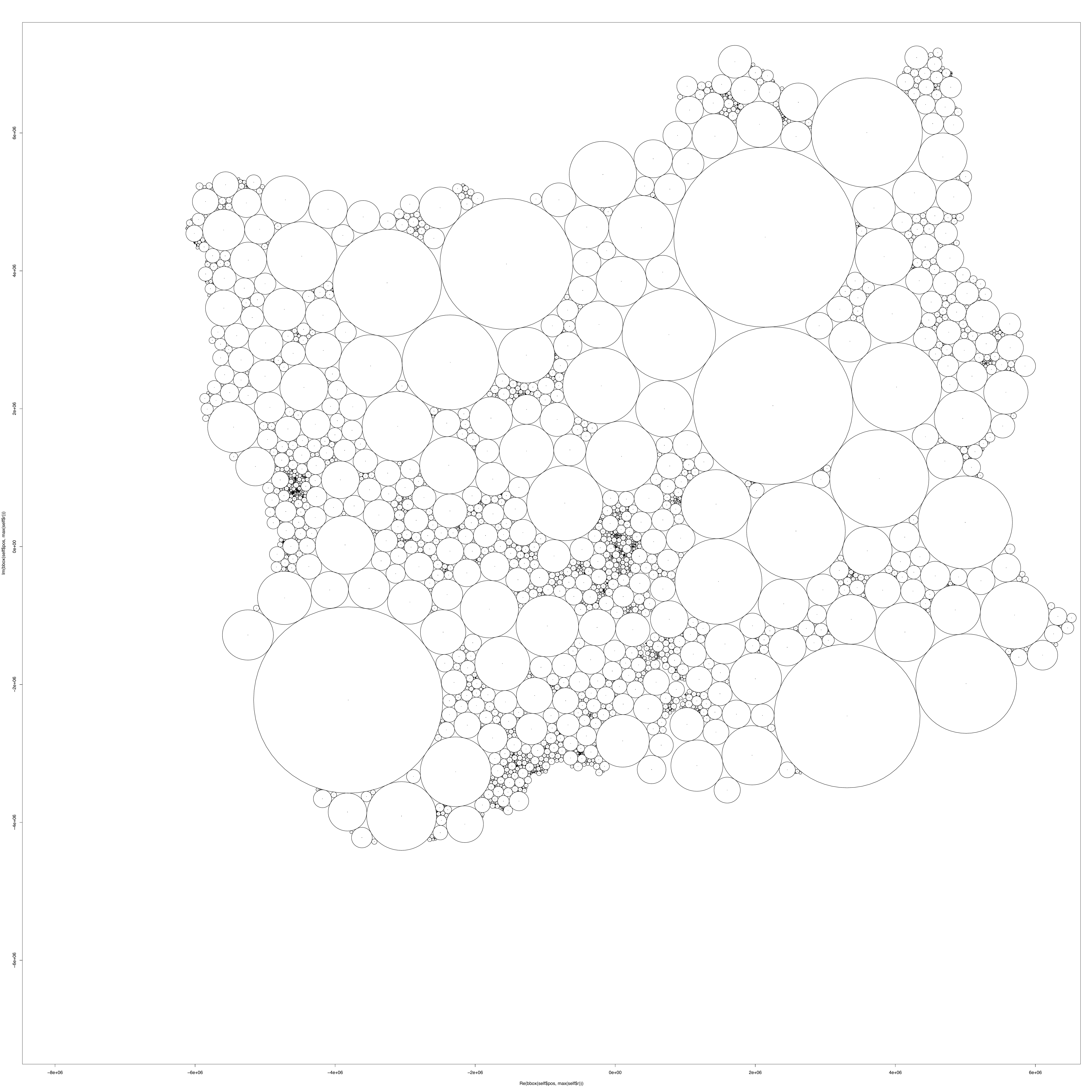}
\end{center}
\caption{This picture is a simulation by Maxim Krikun which represents a circle packing of a uniform triangulation on the disc. }
\label{f.cp}
\end{figure}

If one considers planar maps $\m_n \in \M_n^3$ (i.e. triangulations of the sphere), then by K\"obe Theorem, there is a unique (again up to M\"obius transformations) circle packing in $\S^2$ whose connectivity graph corresponds to $\m_n$.  See figure \ref{f.cp}.

This circle-packing embedding is different from the one given by the conformal structure, yet one might conjecture that if $\m_n$ is sampled uniformly from $\M_n^3$, then the two embeddings should look almost alike with large probability.

\subsubsection{Scaling limit of planar maps embeddings}

Let $\widehat \m_n \subset \S^2$ be a ``natural'' embedding of a random planar map $\m_n\in\M_n^p$ into $\S^2$ ($\m_n$ may be sampled either uniformly as in subsection \ref{ss.PM} or weighted according to some statistical physics model as in subsection \ref{ss.CM}). To each embedding $\widehat \m_n$, the renormalized graph distance ${n^{-a}} d_{gr}$ on the vertices of $\widehat \m_n\subset \S^2$ can be easily extended to a metric $\rho_n$ on  the whole sphere $\S^2$.
 \,If $\m_n$ is sampled uniformly in $\M_n^p$, then we know from Theorem \ref{th.u} that one should choose $a=1/4$. 
If one could prove in this case that as $n\to \infty$, the (random) metric $\rho_n$ on $\S^2$ would converge in law, then it would give a ``good'' answer to Question \ref{NQ} and it would enable us to build a setup for a concrete interpretation of the KPZ formula when $\gamma=\sqrt{8/3}$.

In some sense the idea of Duplantier and Sheffield is to focus on {\bf measures} instead of distances. If $\widehat \m_n\subset \S^2$ is a natural embedding of a planar map $\m_n$, then it is natural to consider the pushforward in $\S^2$ of the Lebesgue measure on $\m_n$, renormalized so that $\m_n$ has unit area (i.e. all faces of $\m_n$ have Lebesgue measure exactly $n^{-1}$). Let us denote by $\mu_n$ this pushforward measure on $\S^2$. $\mu_n$ is thus a random measure on the sphere with $\mu_n(\S^2)=1$. 
As one can see from Figure \ref{f.cp}, one expects that $\mu_n$ should look quite singular with respect to Lebesgue measure on $\S^2$
(or $\D$ if one considers maps with the topology of a disc). Let $f_n$ be the Radon-Nikodym derivative of $\mu_n$ with respect to the Euclidean Lebesgue measure, then we expect $f_n$ to become more and more ``rough'' as $n\to \infty$. 
This brings us to the question. 

\begin{question}\label{RN}
\ni
\bi
\item[(i)] If $\m_n$ is sampled uniformly in $\M_n^p$, is it the case that $\mu_n$ converges as $n\to \infty$ to a random measure on $\S^2$ ?
\item[(ii)] What if  $\m_n\in \M_n^4$ is weighted by $Z_{J_c}(\m_n)$ ?
\ei
\end{question}

Solving this question is a major open problem in the area, but Duplantier and Sheffield made a decisive step in this direction: they managed to {\bf identify} 
an explicit candidate for the scaling limit of $\mu_n$, for which the KPZ equation holds (in the sense of measures). This candidate as we will see in more details below is given by the {\it exponential of a Gaussian Free Field}.

\subsubsection{Setup and statement of the main Theorem}

The common base space will be $[0,1]^2$ and we will consider some deterministic subset $K\subset [0,1]^2$. 
Now, we consider two measures on this space:
\bi
\item[(i)] The Lebesgue measure $\L$ on $[0,1]^2$
\item[(ii)] A random measure $\mu=\mu_\gamma$ which can be formally written as $\frac {d\mu}{d\L} = e^{\gamma\, h}$, where $h$ is an instance of a {\bf Gaussian Free Field}. 
\ei

\begin{figure}[!htp]
\begin{center}
\includegraphics[width=\textwidth]{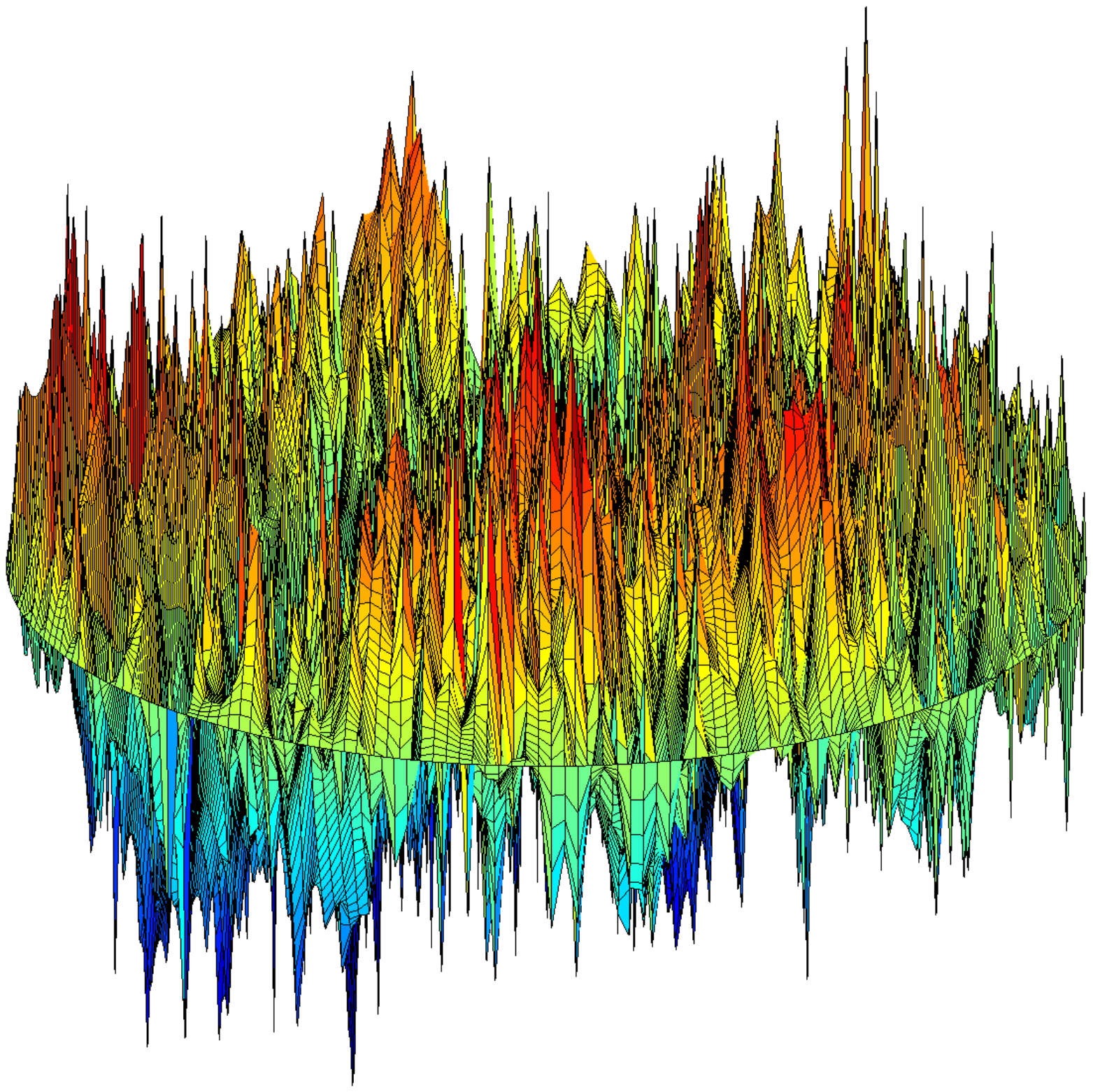}
\end{center}
\caption{An instance of a Gaussian Free Field $h$ in the disc $\D$ with Dirichlet boundary conditions. Picture by  Nam-Gyu Kang. }\label{f.gff}
\end{figure}

The Gaussian Free field (GFF for short) is a certain Gaussian process which will be defined in section \ref{s.gff}. It is ``too rough'' to live in the space of functions on $[0,1]^2$ (say $L^2([0,1]^2)\,$) and  has to be viewed instead as a random Schwartz distribution. See Figure \ref{f.gff} for a
representation of how a regularized GFF on $\D$ looks like. 

\begin{figure}[!h]
\begin{minipage}[t]{0.48\textwidth}
\centering
\includegraphics[width=0.9\textwidth]{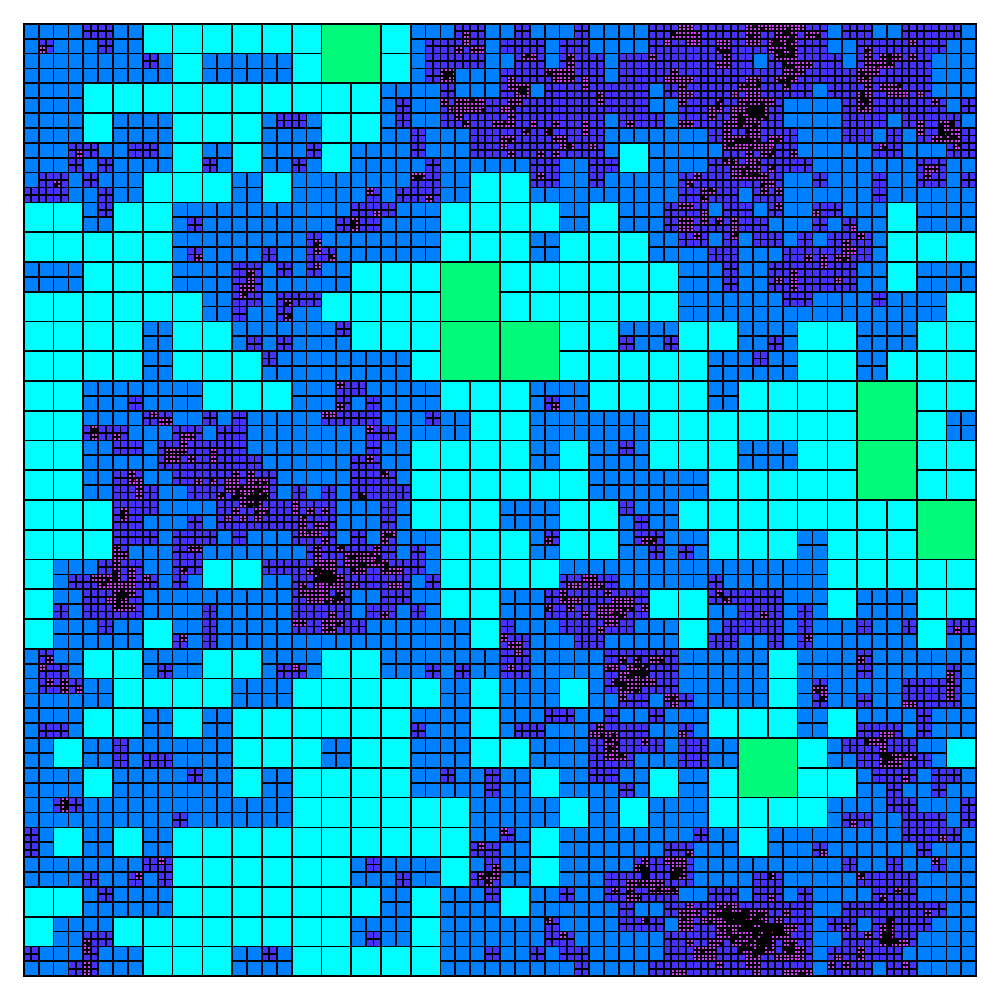}
\end{minipage}
\hspace{0.03\textwidth}
\begin{minipage}[t]{0.48\textwidth}
\centering
\includegraphics[width=0.9\textwidth]{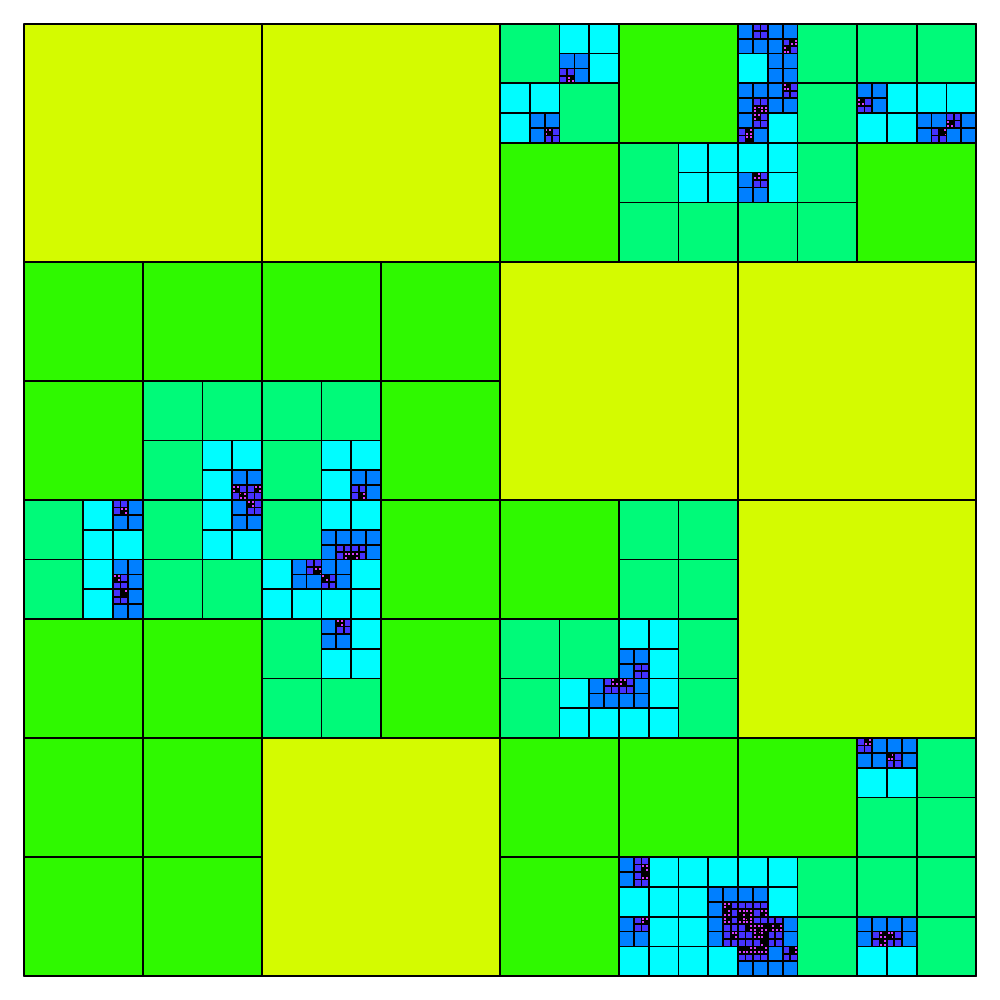}
\end{minipage}
\caption{In each picture, all the plain dyadic squares have about the same quantum area $\delta$. The picture on the left is for $\gamma=3/2$ while the picture on the right which looks much more singular is for $\gamma=5$.}\label{f.QM}
\end{figure}

Since $h$ is not a proper function, $e^{\gamma h}$ is not well-defined a priori. This will be the content of section \ref{s.M} to give a meaning to such measures.
See figure \ref{f.QM} for an illustration of these measures when $\gamma\in \{3/2,5\}$.
These random measures are supposed to model the effect of quantum gravity on our base-space $[0,1]^2$. They are called {\bf Liouville measures}. 

Let us now define Euclidean and quantum scaling exponents in this setting. 
For this, following the notations of \cite{\DS}, we will need the following notion of Euclidean and quantum balls.
\begin{definition}[Euclidean and quantum balls]\label{d.QB}
For any point $z\in [0,1]^2$ and any $\eps,\delta>0$, define 
\bi
\item[(i)] $B_\eps(z)$ to be the Euclidean ball of radius $\eps$ around $z$. 
\item[(ii)] $B^\delta(z)$ to be the ball $B_\tau(x)$ with $\tau:= \sup \{r\geq 0,\, \mu(B_r(z))\leq \delta\}$. $B^\delta(z)$ will be called the {\bf quantum ball} around $z$ of quantum area $\delta$. 
\footnote{Note that since we do not have at our disposal any ``quantum distance'' on $\S^2$, our definition of quantum balls still relies somewhat on the Euclidean one. Nevertheless, it is believed that this slight Euclidean use should ``average out'' as $\delta\to0$. 
}
\ei
\end{definition}

\begin{definition}[scaling exponents]\label{}
Let $K\subset [0,1]^2$ be fixed. 
\begin{enumerate}
\item The {\bf Euclidean scaling exponent} $x=x(K)$ is defined as 
\begin{equation}\label{e.ese}
x=x(K):= \lim_{\eps\to 0} \frac{\log \Pb{B_\eps(z)\cap K \neq \emptyset}}{\log \eps^{2}}\,,
\end{equation}
where the point $z$ is sampled uniformly on $[0,1]^2$.
 
\item The {\bf quantum scaling exponent} $\Delta=\Delta(K)$ is defined as
\begin{equation}\label{e.qse}
\Delta=\Delta(K):= \lim_{\delta\to 0} \frac{\log \Eb{  \mu \big[ B^\delta(z)\cap K \neq \emptyset \big]  } }{\log \delta}\,,
\end{equation}
where $z$ is ``sampled'' according to the a.s. finite measure $\mu=e^{\gamma h}$  and $\E$ averages over the random measure $\mu$.  
\end{enumerate}
As we did before, we implicitly assume here that the limits exist. 
\end{definition}

Note that these definitions are the exact analogs of equations~\eqref{e.ese} and~\eqref{e.qse}. Indeed if one believes in $\mu_n \to \mu$ and if one approximates $K$ be a subset $K_n\subset \widehat \m_n$, then with $\delta=1/n$ one has 
\begin{eqnarray*}\label{}
\Delta=\Delta(K) &\approx& \frac{\log \Eb{  \mu \big[ B^\delta(z)\cap K \neq \emptyset \big]  } }{\log \delta} \\ & \approx &
\frac{\log \Eb{ \mu_n \big[  \text{ the face containing $z\sim \mu_n$ in $\widehat \m_n$ is contained in }K_n \big] }  }{\log 1/n}  \\
& = & \frac{\log \Eb { |K_n| /n } }{\log 1/n}\,,
\end{eqnarray*}
and similarly for $x=x(K)$. We are now in position to state the main theorem of Duplantier and Sheffield.

\begin{theorem}[Duplantier and Sheffield \cite{\DS}]\label{th.main}
Consider the Liouville measure $\mu_\gamma=e^{\gamma h}$ on the unit square $[0,1]^2$ and let $K$ be a (deterministic) subset of $[0,1]^2$, such that the limits in~\eqref{e.ese} and~\eqref{e.qse} exist.
Then if $\gamma\in[0,2)$, 
the quantum scaling exponent $\Delta=\Delta(K)$ for the Liouville measure $\mu_\gamma$ almost surely satisfies to the KPZ formula:
\begin{equation}\label{e.KPZth}
x=\frac{\gamma^2} 4 \Delta^2 + (1-\frac {\gamma^2} 4) \Delta\,.
\end{equation}
\end{theorem}

\begin{remark}\label{}
Note that in the statement of the theorem, by Fubini's theorem, $K$ may also be a random subset of $[0,1]^2$ but in that case, it needs to be chosen independently of $\mu_\gamma$. 
\end{remark}

Their theorem gives a concrete setup in which the KPZ formula holds and enabled Duplantier and Sheffield to 
state the following striking conjecture.

\begin{conjecture}[Duplantier, Sheffield \cite{\DS}]\label{conj}
If $\m_n\in \M_n^p$ are sampled according to a statistical physics model in the $\gamma$-universality class, then 
the pushforward measures $\mu_n$ of any ``natural'' embedding $\widehat \m_n\subset \S^2$ of $\m_n$ weakly converge as $n\to \infty$ towards a random measure which is closely related to the Liouville measure $\mu_\gamma = e^{\gamma h}$, where $h$ is an instance of the Gaussian Free Field on the sphere $\S^2$
\footnote{If one keeps track of the root, asymptotically it will be distributed according to $\mu_\gamma=e^{\gamma h}$.}.
See remark \ref{r.GFFS} for a definition of the GFF on the sphere $\S^2$.
\end{conjecture}
The reason why the limiting measure is not given exactly by the Liouville measure $e^{\gamma h}$ stands from the fact that the measures $\mu_n$ are renormalized to have measure one, while $e^{\gamma h}$ has a random a.s. finite total mass. The actual limit is not given simply by conditioning $e^{\gamma h}$ to have measure one, nor by renormalizing by $\int e^{\gamma h}$, it is slightly more subtle than that.  See section 6 in \cite{\Scott}
for a precise conjecture on the limiting measure.

\medskip

The rest of this survey is divided as follows. In section \ref{s.gff}, we will give a short introduction to Gaussian Free Field. Then, along section \ref{s.M}, we will give a meaning to the Liouville measures $\mu_\gamma=e^{\gamma h}$. In \cite{\DS}, the Liouville measures are defined for all $\gamma\in[0,2)$, we will give here a simplified proof which holds only for $\gamma\in [0,\sqrt{2})$. Finally, in section \ref{s.MP}, we will sketch the main ideas behind the proof of the main theorem. 

\medskip

\ni
{\bf Acknowledgments:}

I wish to thank Vincent Beffara, Nathanael Berestycki, C\'edric Bernardin, Nicolas Curien, Bertrand Duplantier, Gr\'egory Miermont, Rémi Rhodes, Mikael De La Salle, Scott Sheffield, Vincent Vargas and Wendelin Werner for very useful discussions and Bertrand Duplantier,  Nam-Gyu Kang, Maxim Krikun and Scott Sheffield for some of the nice pictures that appear here.

\section{The Gaussian Free Field (GFF)}\label{s.gff}

We will give here a short and by no means self-contained introduction to the Gaussian Free Field. We refer to \cite{\Dub, \GFF} for complete references on this topic. We will try to give a certain flavor of what the GFF is and at the same time, we will introduce some of its key properties which will be needed later.

\subsection{Discrete Gaussian Free Field (DGFF) in the square}
Let us start by a discrete version, known as the Discrete Gaussian Free Field. To simplify, we will consider the case of the square domain $[0,1]^2$. 

\begin{figure}[!htp]
\begin{center}
\includegraphics[width=0.7\textwidth]{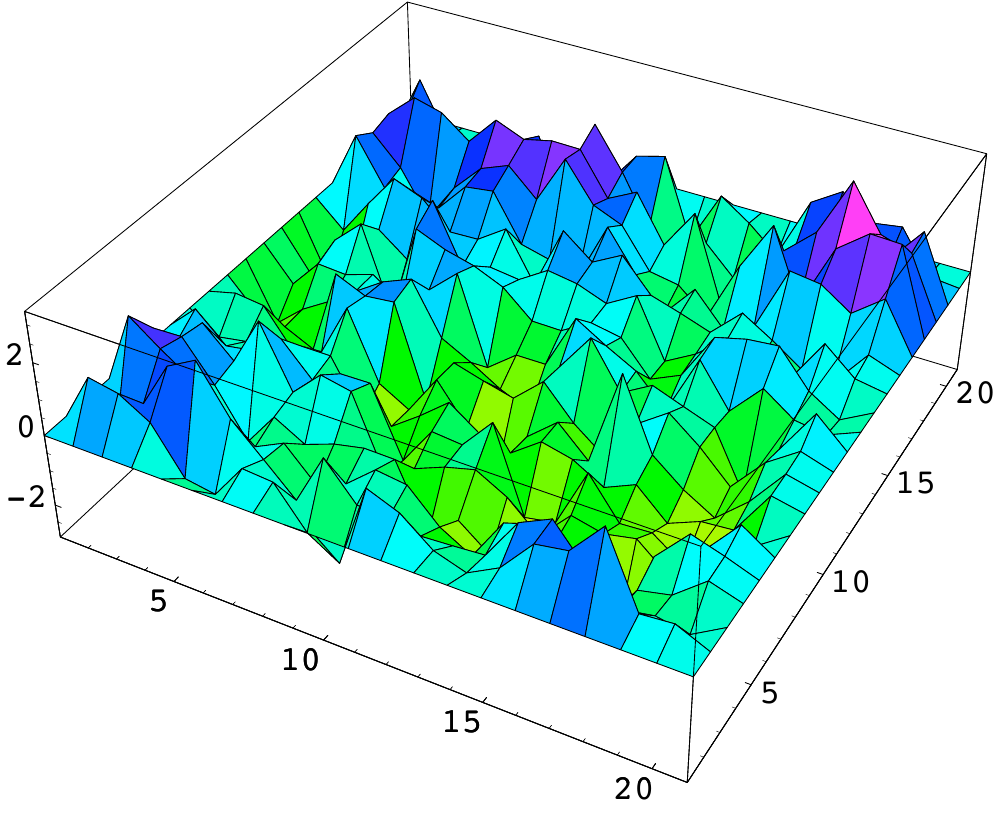}
\end{center}
\caption{A sample of a DGFF on $\Lambda_{21}$. Picture by S. Sheffield.}\label{f.DGFF}
\end{figure}

\begin{definition}\label{d.DGFF}
The DGFF in $\Lambda_N:= \frac 1 N \Z^2 \cap [0,1]^2$ with {\bf Dirichlet boundary conditions} is a probability measure on functions $h_N: \Lambda_N \to \R$ such that 
$h_{ | \,\p \Lambda_N} = 0$ (where $\p \Lambda_N:= \frac 1 N \Z^2 \cap \p [0,1]^2$) and with density
\begin{equation*}\label{}
d\Pb{h_N} := \frac 1 {Z}\,  \exp \Bigl( - \frac 1 2 \sum_{x\sim y} (h_N(x) - h_N(y))^2 \Bigr) \prod dh_N(x)\,,
\end{equation*}
where the sum is over nearest neighbor pairs $x\sim y \in \Lambda_N$, and $Z$ is a renormalizing constant to make it a probability measure.  See figure \ref{f.DGFF} for a sample of a DGFF.
\end{definition}

Due to the Dirichlet boundary conditions, the density is also equal to $Z^{-1} \exp \big( 1/2  \<{h_N, \Delta h_N}\big)$, where $\Delta$ is the discrete Laplacian. Written this way, one sees that DGFF is a {\bf Gaussian random surface} with covariance matrix
given by $\Delta^{-1}$. It is a standard fact that $\Delta^{-1}$ is given by the matrix $[G_N(x,y)]_{x,y\in \Lambda_N}$, where $G_N(x,y)$ is the Green's function for the Random Walk in $\Lambda_N$ killed on $\p \Lambda_N$.  (A good reference for discrete Green's functions is for example the book \cite{\lawler}).
As such, one may give the following equivalent definition of DGFF:

\begin{definition}\label{d.DGFF2}
The DGFF on $\Lambda_N$ with D.b.c is the {\bf centered Gaussian process} $h_N$ indexed by the points $x\in \Lambda_N$ and with covariance structure given by 
\begin{equation*}
\Cov{h_N(x),h_N(y)} = \Eb{h_N(x) h_N(y)} := G_N(x,y)\,.
\end{equation*}
\end{definition}
\begin{remark}\label{r.logN}
Just to give an idea of the amount of fluctuations, ``in the bulk'', say at $x_N:=(N/2,N/2)$, one has 
\begin{equation}\label{e.logN}
\Var{h_n(x_N)}= G_N(x_N,x_N) \asymp \log N\,.
\end{equation}
\end{remark}

\begin{remark}\label{}
If $D\subsetneqq \C$ is any domain of the plane, one can define similarly a DGFF $h_N$ with D.b.c. in the domain $D$ by using the discretization $D_N:= \frac 1 N \Z^2 \cap D$ (and with boundary defined as $\p D_N := \frac 1 N \Z^2 \cap D^c$). 
\end{remark}

It turns out that as $N\to \infty$, $h_N$ converges (in a certain sense to be precised later)
towards a conformally invariant object called {\bf Gaussian Free Field}.

\subsection{A first attempt at defining the Gaussian Free Field}\label{s.informal}
If we are following Definition \ref{d.DGFF2}, it is tempting to define the continuous limit of $h_N$ in the domain $D$ as the centered Gaussian process 
indexed by the points $x\in D$ with the following covariance structure:
\begin{equation}\label{e.informal}
\Cov{h(x),h(y)}:= G_D(x,y)\, \text{ for all }x,y\in D\,,
\end{equation}
where $G_D(x,y)$ is the Green's function of the domain $D$. See the later subsection \ref{s.GF} for a definition (and more) on continuous Green's functions. 
Unfortunately, this definition is ``ill-posed'' since, as suggested by Remark \ref{r.logN}, one would have for any $x\in D$
\begin{equation*}
\Var{h(x)}:=G_D(x,x)=\infty\,.
\end{equation*}

One way to overcome this would be to {\it regularize} $h$ using smooth functions, and this is what we will do when we will introduce the $\eps$-regularization $h_\eps(z)$. Before, let us follow a different approach to define GFF inspired by our initial definition \ref{d.DGFF}.

\subsection{GFF as a Gaussian process indexed by Sobolev space $\Sob^1$}

Following Definition \ref{d.DGFF}, it is natural to look for a probability measure on functions $h:\bar D \to \R$ satisfying $h_{| \p D}\equiv 0$  whose intensity, informally would be given by
\begin{equation}\label{e.inf}
\Pb{h}\propto \exp \Big( - \frac 1 2 \int_D \| \nabla h \|^2 \Big)\,. 
\end{equation}

In order to find a well-defined object corresponding to this informal definition, it is useful to introduce the following Hilbert space:
\begin{definition}[The sobolev space $\Sob^1$]\label{}
Let $\mathcal{C}=\mathcal{C}_D$ be the set of smooth functions $f$ with compact support in $D$.
We define the space $\Sob^1(=\Sob^1_0(D))$ as the closure of $\mathcal{C}$ for the norm 
\begin{equation*}\label{}
\|f\|_\nabla^2:=\frac 1 {2\pi} \int_D \|\nabla f\|^2\,.
\end{equation*}
$\Sob^1$ is a separable Hilbert space for the scalar product
 \begin{equation*}\label{}
\<{f, g}_\nabla:= \frac 1 {2\pi} \int \nabla f\,\nabla g\,.
\end{equation*}
\end{definition}

We are thus trying to define a random function $h\in \Sob^1$ such that \footnote{We do not pay attention to the constant in the exponential through this informal discussion.}
\begin{equation*}
\Pb{h}\propto \exp \Big( - \frac 1 2 \| h \|_\nabla^2 \Big)\,. 
\end{equation*}
 
To gain some intuition, if the Hilbert space $\Sob^1$ happened to be finite dimensional, it would be some $(\R^k, \|\cdot \|_2)$, with $k\geq 1$. 
In that case, if $\mathbf{e}_1,\ldots,\mathbf{e}_k$ is any orthonormal basis of $\R^k$, then a random variable  $x \in \R^k$ with intensity  $\Pb{x}\propto e^{-1/2 \|x\|_2^2}$ can always be written
\begin{equation*}\label{}
x=\sum_{i=1}^k a_i \, \mathbf{e}_i\,,
\end{equation*}
where $(a_i)$ are independent Gaussian random variables. 

Since $\Sob^1$ is a separable Hilbert space, let $(\mathbf{e}_n)_{n\geq 1}$ be an orthonormal basis for $\Sob^1$. By analogy with the finite dimensional case, we want to define the Gaussian Free Field $h$ as 
\begin{equation}\label{e.h}
h:=\sum_{n\geq 1}\, a_n \,\mathbf{e}_n\,,
\end{equation}
where $(a_n)$ are independent Gaussian variables $\sim \mathcal{N}(0,1)$. The difference with the finite-dimensional case is that 
for any $k\geq 1$, the above Gaussian random variable $x$ was almost surely in $\R^k$, while in our present case, it can be shown that almost surely, the above formal series~\eqref{e.h} does not converge in $\Sob^1$. 

In fact, it does not even converge in $L^2(D)$, and as such the Gaussian Free Field $h$ will not be defined as a proper function, but instead as a {\it generalized function} in $\mathcal{D}'$ (i.e. a Schwartz distribution). More precisely,
it can be shown (see \cite{\Dub}) that almost surely, the above sum $h$ converges in the space $\Sob^{-1}$ \footnote{In fact, it turns out that a.s. $h\in \Sob^{-\eps}$ for any $\eps>0$}. 

\begin{definition}\label{d.gff}
From now on, the Gaussian Free Field with Dirichlet b.c. in a domain $D$ will be defined as the random distribution
\begin{equation*}\label{}
h:=\sum_{n\geq 1} a_n \, \mathbf{e}_n   \text{  a.s in }\Sob^{-1}\,,
\end{equation*}
where $(\mathbf{e}_n)_n$ is an orthonormal basis of $\Sob^1(D)$. (The Dirichlet b.c. is hidden in the fact that any $f\in \Sob^1$ satisfies $f_{|\p D}=0$). 
\end{definition}

\begin{example}\label{e.square}
In the case where $D$ is the square $[0,1]^2$, one can write down an explicit basis for $\Sob^1{([0,1]^2)}$, namely 
 for all $j,k \in \N^*$, let
\begin{align}
\mathbf{e}_{j,k}(x,y):= \frac{1}{\sqrt{j^2+k^2}} \, 2\sqrt{2\pi}\, \sin(j \pi x)\, \sin(k \pi y) \,.
\end{align}
It is not hard to check that $(\mathbf{e}_{j,k})_{j,k\in\N^*}$ is indeed an orthonormal basis for $(\Sob^1, \| \cdot \|_\nabla)$.
A Gaussian Free Field in the square $[0,1]^2$ with zero boundary conditions can thus be written as 
\begin{equation*}\label{}
h=\sum_{j,k\in\N^*} \frac{a_{j,k}}{\sqrt{j^2+k^2}} \, 2\sqrt{2\pi}\, \sin(j \pi x)\, \sin(k \pi y) \,,
\end{equation*}
where $(a_{j,k})_{j,k}$ are independent Gaussian variables of variance one and where the convergence for this series holds in the space $\Sob^{-1}$.
\end{example}

\begin{definition}\label{}
If $f = \sum_{n\geq 1} \alpha_n \mathbf{e}_n \in \Sob^1$, then with a slight abuse of notation we will denote by $\<{h,f}_\nabla$ the following quantity 
\begin{equation*}\label{}
\<{h,f}_\nabla:= \sum_{n\geq 1}\, a_n \, \alpha_n\,. 
\end{equation*}
\end{definition}

It is straightforward to check that for any $f\in \Sob^1$, $\<{h,f}_\nabla$ is a Gaussian variable. More precisely, the following 
proposition follows easily from the definition of $\<{h,f}_\nabla$. 
\begin{proposition}\label{pr.GProc}
Let $h$ be a Gaussian Free Field in $D$. Then the process $\big(\<{h,f}_\nabla\big)_{f\in \Sob^1}$ is a centered Gaussian process indexed by the set $\Sob^1$ and with covariance structure
\begin{equation*}\label{}
\Cov{\<{h,f}_\nabla, \<{h,g}_\nabla} = \<{f,g}_\nabla \text{  for any }f,g\in \Sob^1\,.
\end{equation*}
\end{proposition}

\begin{remark}\label{}
In fact, this proposition can serve as another way to introduce the Gaussian Free Field. This is for example the point of view in \cite{\Dub,\GFF}, where they introduce GFF as this Gaussian process indexed by $\Sob^1$. 
\end{remark}

This approach thus gives a good generalization of definition \ref{d.DGFF} to the continuous setting. In fact it can be proved that 
the DGFF $h_N$ weakly converges in the sense of distributions (for example in $\Sob^{-1}$) towards $\lambda_{\Z^2}\,h$, where $\lambda_{\Z^2}$ is a lattice-dependent constant. 

The following proposition relates the Gaussian Free Field we have just defined with definition \ref{d.DGFF2}.
\begin{proposition}\label{pr.GD}
Let $h$ be a Gaussian Free Field in $D$. For any $\rho \in \mathcal{C}_D$ (the smooth functions with compact support in $D$),
we will denote by $\<{h,\rho}$ the distribution $h$ tested against the smooth function $\rho$. Then the process $\big( \<{h,\rho}\big)_{\rho\in \mathcal{C}_D}$ is a centered Gaussian process indexed by $\mathcal{C}_D$ with covariance matrix
\begin{equation*}\label{}
\Cov{\<{h,\rho}, \<{h,\rho'}} = \iint_{D\times D}  \rho(x)\rho'(y)\, G_D(x,y)\,dxdy\,,
\end{equation*}
where $G_D$ is the {\bf Green's function} of domain $D$.
\end{proposition}

\subsection{The Green's function $G_D$ of a domain $D$}\label{s.GF}

\begin{definition}[Green's function in the domain $D$]\label{}
The Green's function of a domain $D\subsetneqq \C$ will be denoted by the function $G_D : D\times D \to \R_+$.
For any $x\in D$, define the function $G^x(y):=G_D(x,y)$.
With such notations, the Green's function $G_D(x,y)$ is characterized by the following properties
\bi
\item[(a)] $\Delta G^x(\cdot) = 0$ on $D\setminus \{x\}$, namely it is harmonic in $D\setminus \{x\}$. 
\item[(b)] $G^x(y)\to 0$ as $y\to \p D$. 
\item[(c)] $G^x(y) \sim \log \frac {1}{|x-y|}$ as $y\to x$. 
\ei
By removing the logarithmic singularity, one can rewrite $G^x(y)$ in the following way $G^x(y)=\log  \frac {1}{|x-y|} + \tilde G^x(y)$, where $\tilde G^x(y)$ is the harmonic extension to $D$ of the function $-\log  \frac {1}{|x-y|}$ on the boundary $\p D$. 
\end{definition}

Here are some well-known properties of Green's functions that we will use.
\begin{proposition}[Properties of the Green's function]\label{}
\ni
\bi
\item[(a)] Conformal invariance: if $\phi : D\to D'$ is a conformal map, then for any $x,y\in D$, 
\[
G_D(x,y)=G_{D'}(\phi(x), \phi(y))\,.
\]
This follows easily from the definition of Green's function.
\item[(b)] Note that $G_\D(0,y)=\log \frac 1 {\|y\|}$.
\item[(c)] $G_D(x,y)=G_D(y,x)$ (this can be seen for example using (a) with (b)).
\item[(d)] For any $x\in D$, $G^x\in \Sob^1$. 
\item[(e)] In the sense of distributions, $\frac {-1}{2 \pi} \Delta[G^x(\cdot)] = \delta_x$, the Dirac point mass at $x$. 
\item[(f)] The above harmonic correction $\tilde G^x(y)$ satisfies 
\begin{equation*}\label{}
\tilde G^x(x)=\log C(x,D)\,,
\end{equation*}
where $C(x,D)$ is the {\bf conformal radius} of $D$ viewed from $x$. If $\phi$ is a conformal map $D \to \D$ with $\phi(x)=0$, then $C(x,D)$ is simply defined as $|\phi'(x)|^{-1}$. This property $(f)$ can be easily checked using (a) and (b). 
\ei
\end{proposition}

Let us now explain how one can recover (at least formally) Proposition \ref{pr.GD} from the above property $(e)$. For any $f\in \mathcal{C}_D$,
one can make sense of $\frac 1 {2\pi} \<{\nabla h, \nabla f}$ where $\nabla h$ is understood in the sense of Schwartz distributions. It is not hard to check that this quantity is exactly $\<{h,f}_\nabla$. Now, since $f$ has compact support, integration by parts implies 
\begin{equation*}\label{}
\<{h,f}_\nabla = \frac 1 {2 \pi} \<{\nabla h, \nabla f} = \frac 1 {2\pi} \<{h, [-\Delta] f}\,.
\end{equation*}
Using this identity with $\rho:= [-\Delta]f \in \mathcal{C}_D$, we find 
\begin{equation*}\label{}
\<{h,\rho} =  \<{h, [-2\pi \Delta^{-1}] \rho}_\nabla\,.
\end{equation*}
This implies Proposition \ref{pr.GD} since 
\begin{align*}\label{}
\Cov{\<{h,\rho}, \<{h,\rho'}} 
&= \Cov{\<{h, [-2\pi \Delta^{-1}] \rho}_\nabla, \<{h,[-2\pi \Delta^{-1}] \rho'}_\nabla} \\
&=\<{[-2\pi \Delta^{-1}] \rho, [-2\pi \Delta^{-1}] \rho'}_\nabla \\
&= \frac 1 {2\pi} \<{\nabla [-2\pi \Delta^{-1}] \rho, \nabla [-2\pi \Delta^{-1}] \rho'} \\
&= \<{\rho, [-2\pi \Delta^{-1}] \rho'} = \iint_{D\times D} \rho(x) \rho(y) G_D(x,y)\, dx dy\,,
\end{align*}
where in the last equality, we used property $(e)$.

Finally, let us mention that using property $(b)$, it is not hard to extract the following striking property for Gaussian Free Field.
\begin{proposition}[conformal invariance]\label{}
Let $\phi:D'\to D$ be a conformal map. If $h$ is GFF in $D$, then $h':=h\circ \phi$ is a Gaussian Free Field in $D'$.
\end{proposition}

\subsection{The $\eps$-regularized GFF $h_\eps$}\label{ss.regul}

The purpose of this subection is to regularize the Gaussian Free Field $h$ in order to obtain a smooth function $h_\eps$.
For this, we will rely on the following $\eps$-regularization of the Green's function. For any $\eps>0$ and any point $x\in D$, let
\begin{equation*}\label{}
G^x_\eps(y):= \log \frac 1 {\eps \vee |x-y|} + \tilde G^x(y)\,.
\end{equation*}
This regularization has the following important property:
\begin{proposition}\label{}
 For any $x\in D$, and any $\eps>0$, $G^x_\eps\in \Sob^1$. Furthermore, in the sense of distributions, one has the following identity 
 \begin{equation}\label{e.gxe}
\frac {-1} {2\pi} \Delta [G^x_\eps(\cdot)] = \nu_{x,\eps}\,,
\end{equation}
where $\nu_{x,\eps}$ denotes the uniform measure on the circle of radius $\eps$ around $x$, $\p B_\eps(x)$. 
\end{proposition}

\begin{remark}\label{}
In fact, with the above definition of $G^x_\eps$, the proposition as stated is not correct when $x$ is close to the boundary ($d(x,\p D)\le \eps$).
To overcome this issue while keeping the same statement for the proposition, the definition of $G^x_\eps$ has to be modified accordingly near the boundary $\p D$. To keep things simple, we choose in this paper to neglect these effects. We refer to \cite{\DS} where this technicality is handled properly. 
\end{remark}

This regularized Green's function enables us to introduce $h_\eps(z)$ the GFF evaluated against $\nu_{x,\eps}$, the uniform measure on $\p B_\eps(x)$. Informally, it corresponds to $h_\eps(z):= \<{h, \nu_{x,\eps}}$. Let us define it as follows 
\begin{definition}\label{}
If $h$ is a sample of a GFF in $D$, then for any $z\in D$, let
\[
h_\eps(z):= \<{h,G^x_\eps}_\nabla\,,
\]
which is well defined since $G^x_\eps\in \Sob^1$. 
\end{definition}
In fact, it corresponds exactly to our informal definition thanks to the following computation:
since $G^x_\eps \in \Sob^1$, one has 
\begin{eqnarray*}\label{}
 h_\eps(z) &=& \<{h,G^x_\eps}_\nabla \\
 &=&\frac 1 {2\pi }\<{h, [-\Delta] G^x_\eps} \\
&=& \<{h,\nu_{x,\eps}} 
\end{eqnarray*}

\subsection{A Brownian motion out of the GFF}
Along this subsection, we will identify a very useful Brownian motion ``within'' the Gaussian Free Field.
We start with the following lemma:
\begin{lemma}\label{l.var}
For any $z\in D$ and any $\eps>0$\footnote{To be self-contained here, one should assume here that $d(z,\p D)\geq \eps$},
one has that 
\begin{equation}
\Var{h_\eps(z)}=\log \frac 1 \eps \, + \log C(z,D)\,,
\end{equation}
where $C(z,D)$ is the conformal radius of $D$ viewed from $z$. 
\end{lemma}

{\bf Proof: }
\begin{align*}
\Var{h_\eps(z)} & = \Var{\<{h,G_\eps^z}_\nabla}  \\
&= \<{G_\eps^z, G^\eps_z}_\nabla \\
&= \<{G^z_\eps, \nu_{z,\eps}} \\
&= \log \frac 1 \eps + \int  \tilde G^z(x) d\nu_{z,\eps}(dx) \\
&=\log \frac 1 \eps + \log C(z,D)\,. \qed
\end{align*}

The following proposition will be of crucial importance in the remaining of this text:
\begin{proposition}\label{pr.BM}
Let $h$ be a GFF with zero-boundary conditions in some domain $D$. For any point $z\in D$, let 
$t_0^z:=\inf\{t\geq 0:\, B_{e^{-t}}(z)\subset D \}$ and let 
\[
Y_t(z):=h_{e^{-t}}(z)\,,
\]
 be the stochastic process defined for any $t\geq t_0^z$. 
(Recall $h_\eps(z)$ denotes the above regularization). 

Then with such notations, the stochastic process  
\begin{equation*}\label{}
\B_t(z):= Y_{t_0^z + t} - Y_{t_0^z},
\end{equation*}
is a standard {\bf Brownian motion}. 
\end{proposition}

{\bf Proof:}
The family of random variables $\{ \B_t(z)\}_{t\geq 0}$ is clearly a Gaussian process. Therefore, it only remains to check that 
for any $0\le s\le t$, $\Cov{\B_s(z),\B_t(z)}= s (=s\wedge t)$. Let $r_0:= e^{-t_0^z}$, $r_1:=e^{-(t_0^z+s)}$ and $r_2:=e^{-(t_0^z+t)}$. Let us first compute
\begin{align*}\label{}
\Cov{h_{r_1}(z), h_{r_2}(z)} &= \Eb{\<{h, G_{r_1}^z}_\nabla \<{h,G_{r_2}^z}_\nabla}\\
&= \<{G_{r_1}^z, G_{r_2}^z}_\nabla \\
&= \<{G_{r_1}^z, \nu_{r_2,z}}\\
&= \log \frac 1 {e^{-t_0^z-s}} + \log C(z,D)\\ &= t_0^z+s + \log C(z,D)\,.
\end{align*}
One can compute in the same way $\Cov{h_{r_0}(z), h_{r_1}(z)}$ and $\Cov{h_{r_0}(z), h_{r_2}(z)}$. This gives us
\begin{align*}
\Cov{\B_s(z), \B_t(z)} &= \Cov{h_{r_1}(z) - h_{r_0}(z),\, h_{r_2}(z) - h_{r_0}(z) } \\
&= t_0^z+s+t_0^z-2*t_0^z \\
&=s \qed 
\end{align*}

Let us conclude this section on the Gaussian Free Field by the following remark.
\begin{remark}[Gaussian Free Field on $\S^2$]\label{r.GFFS}
We have just defined the Gaussian Free Field on a domain $D$ with Dirichlet Boundary conditions. In the same fashion, one can define a Gaussian Free Field on the sphere $\S^2$ (this is needed for example if one wants to make sense of Conjecture \ref{conj}). In this case, the Green's function is given by 
\begin{equation*}
G_D(x,y):= \log \big[ \mathrm{cotan} \frac \theta 2 \big]\, \text{ for all }x,y\in \S^2\,,
\end{equation*}
where $\theta$ denotes the angle between $x$ and $y$. If one wants to define the Gaussian Free Field on $\S^2$ as a Gaussian process similarly as in definition \ref{pr.GProc}, the construction can be done in the same fashion except that in this case, the natural Hilbert space to consider would be the closure for the norm $\|\cdot \|_\nabla$ of the space  $C^0_\infty(\S^2)$ of smooth functions $\phi:\S^2\to \R$ with $\int_{\S^2} \phi(x) dx=0$ (the integral here is with respect to the area measure on $\S^2$).  

\end{remark}

\section{The Liouville measures $e^{\gamma h}$}\label{s.M}

The purpose of this section is to make sense of these Liouville measures $e^{\gamma h}$, which are crucial in the main Theorem \ref{th.main}.
The approach followed in \cite{\DS} is to discretize $e^{\gamma\, h}$ into $e^{\gamma\, h_\eps}$ (where $h_\eps$ is the $\eps$-regularization of the GFF $h$ we have introduced in subsection \ref{ss.regul}) and to then let $\eps\to 0$. As it will become clear below, without renormalization, $e^{\gamma\, h_\eps}$ would diverge in the space of measures. The natural discretization will be the following one:
\begin{definition}\label{d.mue}
For any domain $D$, any $\gamma\geq 0$ and any $\eps>0$, let $\mu_\eps$ be the measure absolutely continuous with respect to Lebesgue measure $\mathcal{L}$ and such that 
\begin{equation}\label{e.mue}
d \mu_\eps(z):=\eps^{\frac {\gamma^2} {2}}\, e^{\gamma\, h_\eps(z)}\,dz\,.
\end{equation}
\end{definition}

Duplantier and Sheffield prove the following proposition in \cite{\DS}:
\begin{proposition}\label{pr.mue}
If $\gamma\in[0,2)$, then for any domain $D$, almost surely as $\eps\searrow 0$ along powers of two, the measures $\mu_\eps$ weakly converge inside $D$ towards a non-degenerate random measure $\mu_\gamma$ which we  will call the {\bf Liouville measure} of parameter $\gamma$. 
The Liouville measure $\mu_\gamma$ is measurable with respect to the Gaussian Free Field $h$ and we will denote it sometimes by $\mu_\gamma=e^{\gamma\, h}$. 
\end{proposition}

\begin{remark}\label{}\ni
\bi
\item[(i)] If $\gamma\in(0,2)$, it can be shown that the Liouville measure $\mu_\gamma$ is a.s. singular with respect to Lebesgue measure.
\item[(ii)] If $\gamma\geq 2$, in some sense things become ``singular''. See for example the work \cite{\GMC} which studies this case. 
\ei
\end{remark}

We give here a new proof 
of this proposition which holds only for the regime $\gamma\in[0,\sqrt{2})$ and furthermore our convergence result will hold only along a certain subsequence $\eps_k \searrow 0$ that we will not make explicit (in \cite{\DS}, it is also along a particular subsequence, but they show that $\eps_k=2^{-k}$ is enough). We believe this proof is interesting in its own since it is slightly different as the one carried in \cite{\DS}, yet it cannot be extended to the range $\gamma\in[\sqrt{2},2)$. See the proof in \cite{\DS} which gives the full range $[0,2)$. 
\medskip

\ni
\underline{\bf Proof in the case $\gamma\in[0,\sqrt{2})$:}

\ni
To simplify, we will restrict ourselves to the case where $D$ is a bounded domain. 
We wish to prove the following proposition:
\begin{proposition}\label{pr.cauchy}
If $\gamma\in[0,\sqrt{2})$, then for any continuous function $\phi:\bar D \to \R$, the sequence of random variables
\begin{equation*}\label{}
\left\{ \mu_\eps(\phi)\right\}_{\eps>0} = \left\{\int_D \phi(z) \mu_\eps(dz) \right\}_{\eps>0}
\end{equation*}
is a Cauchy sequence in $L^2$.
\end{proposition}

Let us first see why this proposition implies Proposition \ref{pr.mue} in the regime $\gamma\in[0,\sqrt{2})$. 
Let $\mathfrak{M}(\bar D)$ be the space of finite positive measures on $\bar D$. It is well-known that this space equipped with the  weak$^*$ topology ( called ``weak convergence of measures'' in Probability theory) is a complete, metrizable, separable space.
Here is an example of a metric on $\mathfrak{M}(\bar D)$ which induces the weak$^*$ topology. Let $(\phi_j)_{j\geq 1}$ be a countable basis of the separable space $(\mathcal{C}^0(\bar D), \|\cdot\|_\infty)$ of continuous functions on $\bar D$ such that $\|\phi_j\|_\infty \le 1$ for all $j$. Then 
\begin{equation*}
d(\eta^1,\eta^2):=\sum_{j\geq 1} \frac {|\eta^1(\phi_j) - \eta^2(\phi_j)|}{2^j}
\end{equation*}
defines a metric on $\mathfrak{M}(\bar D)$ for the weak$^*$ topology. 

Using Proposition \ref{pr.mue}, one can find a subsequence $(\eps_k)_{k\geq 1}$ such that uniformly for all $j\in \{1,...,k\}$, 
\begin{equation*}
\Eb{\Big( \mu_{\eps_{k'}} (\phi_j) - \mu_{\eps_{k''}}(\phi_j)\Big)^2}\, \le 2^{-3k}\,\text{ for all }k''\geq k'\geq k\,.
\end{equation*}

By Markov's inequality, this implies $\Pb{|\mu_{\eps_{k'}} (\phi_j) - \mu_{\eps_{k''}} (\phi_j)|\geq 2^{-k}} \le 2^{-k}$ for all $k\geq 1 $ and $j\le k$. Using Borel-Cantelli lemma, it is an easy exercise to show that this in turn implies that $\mu_{\eps_k}$ is a.s. a Cauchy sequence in $(\mathfrak{M}(\bar D),d)$. Since the later space is complete, we thus obtain an almost sure limit $\mu=\mu_\gamma\in \mathfrak{M}(\bar D)$. Furthermore, since each random measure $\mu_{\eps_k}$ is clearly measurable with respect to the Gaussian Free Field $h$, we obtain that $\mu_\gamma=\lim \mu_{\eps_k}(h)$ is itself measurable with respect to $h$ as a limit in $(\mathfrak{M}(\bar D),d)$ of $h$-measurable measures.

\medskip
\ni
{\bf Proof of Proposition \ref{pr.mue}:}
Let us start with the simpler lemma
\begin{lemma}\label{}
If $\gamma\in [0,\sqrt{2})$, then for any continuous function $\phi: \bar D \to \R$, we have that 
\begin{equation*}\label{}
\Eb{\mu_\eps(\phi)^2} \underset{\eps\to0} {\longrightarrow}  \iint_{D\times D}  \phi(x) \phi(y) \big[C(x,D) C(y,D)\big]^{\gamma^2/2}\, e^{\gamma^2  G_D(x,y)}\,dxdy\,.
\end{equation*}
\end{lemma}

\ni
{\bf Proof of the lemma:}

\begin{align}\label{e.plug}
\Eb{\mu_\eps(\phi)^2}
&= \iint_{D\times D} \phi(x)\phi(y) \,\eps^{\gamma^2}\; \Eb{e^{\gamma h_\eps(x) + \gamma h_\eps(y)}} dx dy \nonumber\\
&= \iint_{D\times D} \phi(x)\phi(y) \,\eps^{\gamma^2}\; e^{\frac {\gamma^2} 2   \Var{h_\eps(x) + h_\eps(y)} } dx dy\,,
\end{align}
where we used the fact that if $X\sim \mathcal{N}(0,\sigma^2)$, then its Laplace transform is given by $\Eb{e^{\gamma X}}= e^{\frac{\gamma^2 \sigma^2} 2}$. 

Now, 
\begin{align*}\label{}
\Var{h_\eps(x) + h_\eps(y)} 
&= \Var{h_\eps(x)} + \Var{h_\eps(y)} + 2\,\Cov{h_\eps(x), h_\eps(y)} \\
&= \Var{h_\eps(x)} + \Var{h_\eps(y)} + 2\,\<{G_\eps^x, \nu_{\eps,y}}\,.
\end{align*}
If $|x-y|>\eps$, we find exactly \footnote{in fact we also need to assume here that $d(x,\p D)\wedge d(y, \p D) \geq \eps$, but we neglect these boundary issues here which are easy to be taken care of}:
\begin{align*}\label{}
\Var{h_\eps(x) + h_\eps(y)} &=  2\log \frac 1 \eps + \log C(x,D)+ \log C(y,D) + 2\,G_D(x,y)\,.
\end{align*}
Let $H(x,y):=\log [C(x,D)C(y,D)] + 2\,G_D(x,y)$, then if $|x-y|\le \eps$, one finds instead the inequality 
\begin{align*}\label{}
\Var{h_\eps(x) + h_\eps(y)} \le 2\log \frac 1 \eps + H(x,y)\,.
\end{align*}
Plugging these into~\eqref{e.plug} gives us 
\begin{align*}\label{}
\Eb{\mu_\eps(\phi)^2} &=  \iint_{|x-y|>\eps}  \phi(x) \phi(y)\big[C(x,D) C(y,D)\big]^{\gamma^2/2} e^{\gamma^2  G_D(x,y)} dxdy \\
& \; + O( \| \phi\|_\infty^2 \, \iint_{|x-y|\le \eps} e^{\frac {\gamma^2}2 H(x,y)} dxdy )
\end{align*}
To conclude the proof of the lemma, one needs to show that the second term goes to zero as $\eps\to 0$ while the first one remains bounded. The key contribution in both cases is what happens when $x\sim y$. In that case, we know that $G_D(x,y)\sim \log \frac 1 {|x-y|}$. This implies that when $x\sim y$, the term $e^{\gamma^2  G_D(x,y)}$ behaves like 
\begin{align*}\label{}
e^{\gamma^2  G_D(x,y)} = \Big| \frac 1 {x-y}\Big|^{\gamma^2+o(1)}\,,
\end{align*}
where $o(1)\to 0$ as $x\to y$. Using the fact that if $\alpha<2$, then $\iint_{D\times D} \Big| \frac 1 {x-y}\Big|^{\alpha} < \infty$, it is an easy exercise to conclude the proof of the Lemma. \qed

Now, let us prove  
Proposition \ref{pr.cauchy}, i.e. that $\{ \mu_\eps(\phi)\}_{\eps>0}$ is a Cauchy sequence for any $\phi$ continuous on $\bar D$. For this, let us estimate for $0<\eta <\eps$:
\begin{align*}
\Eb{\big( \mu_\eps(\phi) - \mu_\eta(\phi) \big)^2} 
&= \Eb{\mu_\eps(\phi)^2} + \Eb{\mu_\eta(\phi)^2}  \\
&\;  - 2 \eps^\fgam \eta^\fgam \iint_{D\times D} \phi(x)\phi(y) 
e^{\frac {\gamma^2} 2   \Var{h_\eps(x) + h_\eta(y)} } dx dy\,.
\end{align*}

Similarly as in the proof of the Lemma, we find that 
\begin{align*}\label{}
\Var{h_\eps(x)+h_\eta(y)}&\le \log \frac 1 \eps + \log \frac 1 \eta + \log(C(x,D)C(y,D))+2\,G_D(x,y)\,, 
\end{align*}
with equality if and only if $|x-y|\geq \eps\vee \eta$ (and $d(x,\p D)\wedge d(y, \p D) \geq \eps$ as well). 
In particular, in the same fashion as above, if $\gamma$ is chosen so that $\gamma<\sqrt{2}$, this implies that 
\begin{align*}
\eps^\fgam \eta^\fgam \iint_{D\times D} \phi(x)\phi(y) 
e^{\frac {\gamma^2} 2   \Var{h_\eps(x) + h_\eta(y)} } dx dy &\\
& \hskip -5 cm \underset{0<\eta<\eps\to0}{\longrightarrow}    \iint_{D\times D}  \phi(x) \phi(y) \big[C(x,D) C(y,D)\big]^{\gamma^2/2}  e^{\gamma^2  G_D(x,y)}\,dxdy\,,
\end{align*}
which thus implies that $\mu_\eps(\phi)$ is indeed a Cauchy sequence in $L^2$. \QED


\begin{remark}\label{}
Another natural approach would be to consider $\eps \mapsto \mu_\eps(\phi)$ as a stochastic process in $\eps\searrow 0$. Since it can be written as 
\begin{align*}\label{}
\mu_\eps(\phi) & = \int_D \phi(x) e^{\gamma h_\eps(x) + \fgam \log \eps} \, dx\,, \\
&= \int_D \phi(x) C(x,D)^{-\fgam} e^{\gamma h_\eps(x) - \fgam \Var{h_\eps(x)}}\,,
\end{align*}
and since for each fixed $x$, $\eps\mapsto  e^{\gamma_\eps(x) - \fgam \Var{h_\eps(x)}}$ is a positive {\bf martingale}, one might be tempted to prove the a.s. convergence of $\mu_\eps(\phi)$, when $\phi\geq 0$, by showing that it is a positive martingale (furthermore the above $L^2$ bounds when $\gamma<\sqrt{2}$ would imply its uniform integrability). Unfortunately, this is not the case. In \cite{\DS}, the authors manage nevertheless to rely on such an approach by looking at a different way to regularize $h$, namely $h_n:= \sum_{i\le n} \<{h,f_i}_\nabla f_i$, where $\{f_i\}$ is some orthonormal basis of $\Sob^1$. In that case $\mu_n(\phi)$ (defined accordingly, see \cite{\DS}) is in this case a positive martingale. 
\end{remark}

\begin{remark}\label{}
Finally, one should point out that such measures had already been constructed within the theory of {\bf Gaussian Multiplicative Chaos}
initially developed by Kahane in \cite{\Kahane}. See \cite{\GCrevisited} for a more general construction of such measures.
\end{remark}



%



\section{Ideas behind the proof of the main Theorem}\label{s.MP}
In this section, we wish to explain where the KPZ formula comes from by giving some of the ideas behind the proof of Theorem \ref{th.main}.

\subsection{Setup}

Let us fix some (deterministic) $K\subset [0,1]^2$ and some parameter $\gamma \in[0,2)$. 
Assuming that the limit exists in equation~\eqref{e.qse}, our goal is to express the quantum scaling exponent $\Delta=\Delta(K)$ as a function of the Euclidean scaling exponent $x=x(K)$. 
Recall that 
\begin{equation*}
\Delta=\Delta(K):= \lim_{\delta\to 0} \frac{\log \Eb{  \mu \big[ B^\delta(z)\cap K \neq \emptyset \big]  } }{\log \delta}\,,
\end{equation*}
where $z$ is ``sampled'' according to the a.s. finite measure $\mu=e^{\gamma h}$ and $\E$ averages over the random measure $\mu$.  
In this setting, we first sample $\mu_\gamma = e^{\gamma \,h}\sim \E$ and conditioned on the measure $\mu_\gamma$, 
we sample $z\sim \mu_\gamma$. For the proof of the main theorem, it will be useful to invert this procedure, i.e. to first sample $z$ 
(according to the correct marginal measure) and then to sample $\mu_\gamma$ conditioned on the value of $z$. 
The coming subsection introduces the right tool for this.

\subsection{Rooted Liouville measures}
As we have just explained, the definition of the exponent $\Delta$ involves the coupling $(\mu_\gamma,z)$, where $z$ is sampled according to the first coordinate $\mu_\gamma$ (which in general is not a probability measure). The law of this coupling can be written $dh\times d\mu_\gamma(z)$. In order to invert this sampling procedure, we will introduce the following probability measure which can be viewed as an $\eps$-regularization of the 
above coupling.
\begin{definition}[Rooted Liouville measure]\label{}
For any $\gamma\in[0,2)$ and any $\eps>0$, let 
\begin{equation*}\label{}
\Theta_\eps:=\frac 1 {Z_\eps} e^{\gamma h_\eps(z)} \, dhdz\,, 
\end{equation*}
be the probability measure on $\Sob^{-1}\times [0,1]^2$, where $Z_\eps$ is a renormalizing constant chosen so that $\Theta_\eps$ is a probability measure. 
\end{definition}

This regularized coupling enables us to make sense of the reversed sampling procedure. First of all, we need to compute the marginal distribution of $\Theta_\eps$ on $z$. It is simply given by $\rho_\eps(z)=Z_\eps^{-1} \E_h\big[ e^{\gamma h_\eps(z)}\big]$ which is 
explicit since using Lemma \ref{l.var} one has that it is proportional to $C(z,D)^{\gamma^2/2}$ when $d(z, \p D)\geq \eps$. In particular, as $\eps\to 0$, the density $\rho_\eps(z)$ converges towards a limiting $\rho(z)\propto C(z,D)^{\gamma^2/2}$. Now, conditioned on $z$, the marginal on $h$ is given by the Gaussian Free Field $h$ weighted by $e^{\gamma\, h_\eps(z)}$ i.e. by $\frac {e^{\gamma\, h_\eps(z)} dh}{\int e^{\gamma\, h_\eps(z)} dh}$. Note that this step would not make any sense without our $\eps$-regularization since the scalar quantity ``$e^{\gamma\, h(z)}$'' is not defined at point $z$. 
It is a standard fact about Gaussian processes (Cameron-Martin theorem), that if $h^z$ is sampled according to $\frac {e^{\gamma\, h_\eps(z)} dh}{\int e^{\gamma\, h_\eps(z)} dh}$, then $h^z \overset{(d)}{=}h+\gamma\, G^z_\eps$. To see this at least heuristically, note that the process $h+\gamma\, G^z_\eps$ is a deterministic translation of the Gaussian process $h$ and it is easy to check that in the finite dimensional case, if $X$ is a standard Gaussian vector in $\R^n$, then for any fixed $u\in \R^n$, the law of $\tilde X= X+u$ is the same as the law of $X$ weighted by $e^{\<{X,u}}$. In the case of the Gaussian Free Field, the scalar product is $\<{\cdot,\cdot}_\nabla$ and we indeed have by definition that $\<{h,\gamma G^z_\eps}_\nabla= \gamma\, h_\eps(z)$.

In particular, the pair $(z,h^z)\sim \Theta_\eps$ can also be sampled as follows: first sample $z$ according to the above marginal distribution $\rho_\eps$(whose density away from $\p D$ is proportional to $C(z,D)^{\gamma^2/2}$), and then conditioned on $z$, let $h^z:=  h+ \gamma\, G_\eps^z$, where $h$ is an independent Gaussian Free Field in $[0,1]^2$. This way, we see that the measures $\Theta_\eps$  converge towards a limiting measure $\Theta$ as $\eps\to 0$, for which the GFF ($h^z$) conditioned on the first component $z\sim \rho$ is sampled according to $h^z:=  h +\gamma\, G^z$.

Using this coupling we have that 
\begin{align}\label{e.last}
\Eb{  \mu \big[ B^\delta(z)\cap K \neq \emptyset \big]  }
&\asymp \Theta\Big[ B^\delta(z)\cap K \neq \emptyset \Big]  \nonumber \\
& \asymp \int_{[0,1]^2} \rho(z)dz \; \P_{h} \big[ B^\delta_{h^z= h+ \gamma\, G^z} (z) \cap K\neq \emptyset \big]\,,
\end{align}
where $B^\delta_{h^z= h+ \gamma\, G^z}(z)$ denotes the quantum ball of quantum area $\delta$ around $z$ in the sense of Definition \ref{d.QB} but with a field sampled according to $h^z:=h+\gamma\, G^z$.

\subsection{An estimate about quantum balls around the root}
For any $z\in [0,1]^2$, in the same fashion as in Proposition \ref{pr.mue}, it can be shown that if $\gamma<2$ and if  $h^z:=h+\gamma\, G^z_\eps$, 
then one can make sense of the Liouville measure rooted at $z$, $\mu_\gamma^z:=e^{\gamma\, h^z}$. 

In \cite{\DS}, the following property is shown.
\begin{proposition}\label{pr.approx}
Let $z\in [0,1]^2$ and let $h^z:= h + \gamma\, G^z$, then almost surely as the Euclidean radius $r \searrow 0$,
\begin{equation*}
\mu_\gamma^z(B_r(z)) \sim c\,  r^{\gamma Q}\, e^{\gamma \, h_r^z (z)}\,.
\end{equation*}
where $Q=Q_\gamma:=2/\gamma+\gamma/2>2$  and where $h_r^z:= h_r+ \gamma\, G_r^z$. ($c=c_\gamma$ is some explicit constant).
See subsection 4.1 in \cite{\DS}.
\end{proposition}

We will not prove this proposition, but instead we will convince ourselves through the computation of the expectation of $\mu_\gamma^z(B_r(z))$ that one can indeed expect such a behavior:

\begin{lemma}\label{}
\begin{equation*}
\Eb{\mu_\gamma^z(B_r(z))} \sim c\, r^{\gamma Q}\,\Eb{e^{\gamma\, h_r^z(z)}}\,,
\end{equation*}
as $r\to 0$ for a certain constant $c=c_\gamma$. 
\end{lemma}

\ni
{\bf Proof:}
First of all, using Lemma \ref{l.var}, one has (if $d(z, \p D)\geq r$):
\begin{align*}
\Eb{e^{\gamma\, h_r^z(z)}} & =  \Eb{e^{\gamma\, h_r(z)}}  e^{\gamma^2 G_r^z(z)}  \\
&= C(z,D)^\fgam r^{-\fgam}\, C(z,D)^{\gamma^2} r^{-\gamma^2} \\
&= C(z,D)^\frac {3 \gamma^2} 2  r^{-3\gamma^2 /2}\,.
\end{align*}

In order to compute $\Eb{\mu^z_\gamma(B_r(z))}$, 
let us approximate $h^z=h+\gamma\,G^z$ into $h^z_\eps:=h_\eps +\gamma\, G^z$:
\begin{align*}
\Eb{\mu^z_\gamma(B_r(z))} 
&=\lim_{\eps\to 0} \, \eps^{\fgam}\, \int_{B_r} \Eb{e^{\gamma h_\eps(x)+ \gamma^2\, G^z(x) }}\, dx \\
&=\int_{B_r}  e^{\fgam \tilde G^z(x)}\, e^{\gamma^2\, G^z(x)} dx   \text{ using Lemma \ref{l.var}}\\
&=\int_{B_r} e^{\fgam \tilde G^z(x)}\, e^{\gamma^2 \log |x-z|^{-1} + \gamma^2 \tilde G^z(x)} \\
&=\int_{B_r} \exp{\big[\frac {3\gamma^2} 2  \tilde G^z(z)+o(1)\big]}\, \frac {1} {|x-z|^{\gamma^2}} \, dx\,,
\end{align*}
as $r\to 0$, since $x\mapsto \tilde G^z(x)$ is continuous. Therefore as $r\to 0$:
\begin{align*}
\Eb{\mu^z_\gamma(B_r(z))} 
&\sim C(z,D)^{3\gamma^2/2} \int_{B_r}   \frac {1} {|x-z|^{\gamma^2}} dx \\
&\sim C(z,D)^{3\gamma^2/2} \int_{u=0}^r 2\pi \,  u^{1-\gamma^2} du  \\
&  \sim c\,  C(z,D)^{3\gamma^2/2}  r^{2-\gamma^2} \\
& \sim c\, r^{\gamma \, Q} \Eb{e^{\gamma \, h}} \QED
\end{align*}

\begin{remark}\label{}
This first moment computation indeed provides some supporting evidence for Proposition \ref{pr.approx}. Yet, such a comparison of first moments is not so natural after all since, as it was pointed out to us by Nicolas Curien, the expected quantum area of $B_r(z)$ diverges as $r\to 0$ when $\gamma\in (\sqrt{2},2]$ (in the above displayed equation, $2-\gamma^2<0$ when $\gamma>\sqrt{2}$). This counter-intuitive phenomenon  is due to the fact that for any $\gamma>0$, the main contribution in $\Eb{e^{\gamma\, h_r^z(z)}}$ does not come from typical properties of $h_\eps(z)$ but follows instead from large deviations events for $h_\eps(z)$. This is why first moments computations are not suitable for studying ``typical'' behavior as one is interested in in the statement of Proposition \ref{pr.approx}.
See subsection 4.1 in \cite{\DS} for a proof of this ``law of large numbers'' type of behavior. 
\end{remark}

Recall that the quantum ball $B^\delta(r)$ for a field $h^z$ is defined as 
$B_\tau(z)$ with $\tau:= \sup \{r\geq 0,\, \mu_\gamma^z(B_r(z))\leq \delta\}$.
The content of Proposition \ref{pr.approx} tells us that $B^\delta(z)$ should be very well approximated by the ball $\tilde B^\delta(r)$ defined as 
$B_{\tilde \tau}(z)$ with 
\begin{equation*}\label{}
\tilde \tau:= \sup  \left\{ r>0\,\text{ s.t. } c\, r^{\gamma Q}\, e^{\gamma\, h_r^z(z)} \le \delta \right\}\,.
\end{equation*}
Plugging this into~\eqref{e.last} gives us 
\begin{align*}\label{}
\Eb{  \mu \big[ B^\delta(z)\cap K \neq \emptyset \big]  }
&\asymp \int_{[0,1]^2} \rho(z) dz\, \P_h \big[  B_{h^z=h+\gamma G^z} ^\delta(z)\cap K \neq \emptyset \big] \\
&\asymp   \int_{[0,1]^2} dz \, \E_h \Big[  1_{ B^\delta_{h+\gamma G^z}(z) \cap K \neq \emptyset}   \Big] \\
&\approx  \int_{[0,1]^2} dz \, \E_h \Big[  1_{ B_{\tilde \tau}(z) \cap K \neq \emptyset}    \Big]  \\
&\asymp \int_{[0,1]^2}  dz \, \E_h \Big[  \Eb{  1_{ B_{\tilde \tau}(z) \cap K \neq \emptyset}   \md \tilde \tau }\Big]\,.
\end{align*}
As we will see below, it is not difficult to show that the law of the random radius $\tilde \tau=\tilde\tau(z)$ at point $z$ depends very little on the point $z\in[0,1]^2$ \footnote{we are neglecting boundary issues here}.
In particular, if $\tilde \P$ denotes this common law for $\tilde \tau$, we have that
\begin{align*}\label{}
\Eb{  \mu \big[ B^\delta(z)\cap K \neq \emptyset \big]  }
&\approx \int_{\R_+} d\tilde \P(\bar r)   \int_{[0,1]^2} dz   1_{B_{\bar r}(z) \cap K \neq \emptyset} \\
&\approx \int_{\R_+} d\tilde \P(\bar r)  \;  {\bar r}^{2\,x(K)} \text{  by definition of $x=x(K)$}\\
&\approx \tilde \E \big[ {\tilde \tau}^{2\,x(K)}\big]\,.
\end{align*}
Therefore it only remains to understand the law of $\tilde \tau$ (in some sense uniformly in the root $z\in[0,1]^2$).
This will be done by identifying a drifted Brownian motion within $h^z:=h+\gamma \, G^z$.

\subsection{Reduction to a large deviation question on Brownian motion}

Let us fix some $z\in (0,1)^2$ and let $t_0\geq 0$ so that $B_{e^{-t_0}}(z) \subset (0,1)^2$.  Then similarly as in Proposition \ref{pr.BM},
if 
\begin{align*}\label{}
W_t(z)&:= h_{e^{-t-t_0}}^z - h_{e^{-t_0}}^z \\
&= h_{e^{-t-t_0}}(z) - h_{e^{-t_0}}(z) + \gamma\, \big[ G^z_{e^{-t-t_0}}(z) - G^z_{e^{-t_0}}(z) \big] \\
&= h_{e^{-t-t_0}}(z) - h_{e^{-t_0}}(z) + \gamma \, t \,,
\end{align*}
then $(W_t(z))_{t\geq 0}$ is a Brownian motion with drift $\gamma$. I.e $W_t \overset{(d)}= B_t+ \gamma\, t$. 
$\tilde \tau$ can be defined using this Brownian motion: indeed recall
\begin{align*}\label{}
\tilde \tau &:= \sup  \left\{ r>0\,\text{ s.t. } c\, r^{\gamma Q}\, e^{\gamma\, h_r^z(z)} \le \delta \right\} \\
&= \exp  -\Big[ \inf \left\{ t>0 \,\text{ s.t. }  c\, (e^{-t})^{\gamma Q}\, e^{\gamma\, h_{e^{-t}}^z(z)} \le \delta   \right\} \Big] \\
&\approx \exp -\Big[ \inf \left\{ t>0 \,\text{ s.t. }  (e^{-t})^{\gamma Q}\, e^{\gamma\, B_t  + \gamma^2\, t } \le \delta   \right\} \Big]  \\
& = \exp -\Big[ \inf \left\{ t>0 \,\text{ s.t. }  B_t  + (\gamma -Q)\, t  \le \frac {\log \delta} \gamma   \right\} \Big] \\
&=  \exp -\Big[ \inf \left\{ t>0 \,\text{ s.t. }  \bar B_t  + a_\gamma \, t  \geq \frac {\log \frac 1\delta} \gamma   \right\} \Big]\,,
\end{align*}
where $\bar B_t:=(-B_t)$ is a standard Brownian motion and where $a_\gamma:=Q_\gamma -\gamma = 2/\gamma-\gamma/2$, which is positive when $\gamma<2$. 
Let $T=T_\delta$ be the stopping time for the drifted Brownian motion $\bar B_t + a_\gamma \,t$ stopped the first time it reaches level $\gamma^{-1} \log \frac 1 \delta  $. Since $\tilde \tau \approx e^{-T_\delta}$, summarizing the above discussion, we obtain that $\Delta=\Delta(K)$ should be given by 
\begin{equation*}\label{}
\Delta = \lim_{\delta\to 0} \frac{\log \Eb{e^{-2 x\, T_\delta}}}{\log \delta}\,.
\end{equation*}

It remains to compute the quantity $\Eb{e^{-2 x\, T_\delta}}$. This is a classical computation for Brownian motion and it works as follows. 
Consider for any $\beta\geq 0$ the process:
\begin{equation*}\label{}
t\mapsto \exp \big( \beta B_t - \frac{\beta^2} 2 t\big)\,,
\end{equation*}
which is a martingale. Using the optional stopping theorem for the stopping time $T_\delta$, we get for any $\beta\geq 0$,
\begin{align*}\label{}
\Eb{\exp \big( - \beta\,a_\gamma\, T_\delta  + \beta/\gamma \log \frac 1 \delta  - \beta^2\, T_\delta /2 \big)} = 1\,,
\end{align*}
which in turn gives 
\begin{equation*}\label{}
\Eb{e^{-2x T_\delta}} = \delta^{\beta/\gamma} \,,
\end{equation*}
if $\beta=\beta_\gamma$ is chosen so that $2x = \beta\, a_\gamma + \beta^2/2$. Since $\Delta=\Delta(K)$ is given ultimately by $\beta_\gamma/\gamma$, this indeed gives us a quadratic relation between $x(K)$ and $\Delta(K)$. One can check that this quadratic relation is the KPZ formula~\eqref{e.KPZth}. \QED

\subsection{Other proofs of a KPZ formula in the literature}

Finally, let us mention that after Duplantier and Sheffield announced their result, other proofs of KPZ formulas have been proved in slightly different settings:
\bi
\item Benjamini and Schramm obtained in \cite{\BS} a simple and enlightening proof of a KPZ formula for multiplicative dyadic cascades in dimension one. 
The advantage of their proof is that it gives a quadratic relation between actual Hausdorff dimensions as opposed to ``expected box-counting dimensions'' in \cite{\DS} Unfortunately, their argument is inherently one-dimensional and if one would extend their argument to higher dimensions, it would no longer deal with proper Hausdorff dimensions. 

\item Rhodes and Vargas proved in \cite{\RV} a KPZ formula in the general setting of Gaussian multiplicative Chaos. Their proof enables to deal with ``stationary'' measures (as opposed to \cite{\BS} which relies on a discrete dyadic division). 
The difference with \cite{\DS} is that their KPZ formula holds for a different notion of dimension (or rather scaling exponent) 
as the one considered in \cite{\DS}. In that sense their work is complementary to the work \cite{\DS}. More precisely, in rough terms, 
if $K\subset [0,1]^2$ and if $\mu$ denotes a measure (for example the Liouville measure), then their notion of scaling exponent is defined as 
\begin{equation*}
\tilde \Delta(K):= \inf \left\{ s\in(0,1],\, \text{s.t.}   \inf_{\text{coverings } K\subset \cup B(x_i,r_i)} \{ \sum_{i} \mu(B(x_i,r_i))^s \} =0 \right\}\,,
\end{equation*}
where the balls $B(x_i,r_i)$ are Euclidean balls of radii $r_i$. 
This notion is very different from the expected box-counting dimension considered in \cite{\DS}. 
Note that in both works \cite{\DS} and \cite{\RV}, the notions of scaling exponents $\Delta, \tilde \Delta$ still rely somewhat on the Euclidean metric. It seems one is still far from a ``true'' KPZ correspondence between Euclidean and quantum {\it metrics}. 
\ei

%






\bibliography{bourbaki}

\begin{thebibliography}{LGM11b}

\bibitem[AC]{AngelCurien}
Omer Angel and Nicolas Curien.
\newblock {P}ercolations on infinite random maps.
\newblock in preparation.

\bibitem[Ang03]{MR2024412}
Omer Angel.
\newblock Growth and percolation on the uniform infinite planar triangulation.
\newblock {\em Geom. Funct. Anal.}, 13(5):935--974, 2003.

\bibitem[Bet11]{bettinelli}
J.~Bettinelli.
\newblock {S}caling limit of random planar quadrangulations with a boundary.
\newblock ar{X}iv:1111.7227, 2011.

\bibitem[BJRV12]{GMC}
J.~Barral, X.~Jin, R.~Rhodes, and V.~Vargas.
\newblock {G}aussian multiplicative chaos and {KPZ} duality.
\newblock ar{X}iv:1202.5296, 2012.


\bibitem[BC12]{BenjaCurien}
Itai Benjamini and Nicolas Curien.
\newblock {S}imple random walk on the uniform infinite planar quadrangulation:
  {S}ubdiffusivity via pioneer points.
\newblock ar{X}iv:1202.5454, 2012.


\bibitem[BS09]{MR2506765}
Itai Benjamini and Oded Schramm.
\newblock K{PZ} in one dimensional random geometry of multiplicative cascades.
\newblock {\em Comm. Math. Phys.}, 289(2):653--662, 2009.


\bibitem[Dav88]{Dav88}
François David.
\newblock Conformal Field theories coupled to 2-{D} gravity in the conformal gauge. 
\newblock {\em Mod. Phys. Lett. A}, 3(17): 1651--1656, 1988.

\bibitem[DKH88]{DK88}
B.~Duplantier and Kwon K.-H.
\newblock Conformal invariance and intersection of random walks.
\newblock {\em Phys. Rev. Lett.}, 61:2514--2517, 1988.

\bibitem[DS11]{MR2819163}
Bertrand Duplantier and Scott Sheffield.
\newblock Liouville quantum gravity and {KPZ}.
\newblock {\em Invent. Math.}, 185(2):333--393, 2011.

\bibitem[Dub09]{MR2525778}
Julien Dub{\'e}dat.
\newblock S{LE} and the free field: partition functions and couplings.
\newblock {\em J. Amer. Math. Soc.}, 22(4):995--1054, 2009.

\bibitem[Dup98]{MR1666816}
Bertrand Duplantier.
\newblock Random walks and quantum gravity in two dimensions.
\newblock {\em Phys. Rev. Lett.}, 81(25):5489--5492, 1998.

\bibitem[Dup04]{MR2112128}
Bertrand Duplantier.
\newblock Conformal fractal geometry \& boundary quantum gravity.
\newblock In {\em Fractal geometry and applications: a jubilee of {B}eno\^\i t
  {M}andelbrot, {P}art 2}, volume~72 of {\em Proc. Sympos. Pure Math.}, pages
  365--482. Amer. Math. Soc., Providence, RI, 2004.

\bibitem[Kah85]{MR829798}
Jean-Pierre Kahane.
\newblock Sur le chaos multiplicatif.
\newblock {\em Ann. Sci. Math. Qu\'ebec}, 9(2):105--150, 1985.

\bibitem[Kaz86]{MR871244}
V.~A. Kazakov.
\newblock Ising model on a dynamical planar random lattice: exact solution.
\newblock {\em Phys. Lett. A}, 119(3):140--144, 1986.

\bibitem[KPZ88]{MR947880}
V.~G. Knizhnik, A.~M. Polyakov, and A.~B. Zamolodchikov.
\newblock Fractal structure of {$2$}{D}-quantum gravity.
\newblock {\em Modern Phys. Lett. A}, 3(8):819--826, 1988.

\bibitem[LG07]{MR2336042}
Jean-Fran{\c{c}}ois Le~Gall.
\newblock The topological structure of scaling limits of large planar maps.
\newblock {\em Invent. Math.}, 169(3):621--670, 2007.

\bibitem[LG11]{LGu}
Jean-Fran{\c{c}}ois Le~Gall.
\newblock {U}niqueness and universality of the {B}rownian map.
\newblock ar{X}iv:1105.4842, 2011.

\bibitem[LGM11a]{MR2778796}
Jean-Fran{\c{c}}ois Le~Gall and Gr{\'e}gory Miermont.
\newblock Scaling limits of random planar maps with large faces.
\newblock {\em Ann. Probab.}, 39(1):1--69, 2011.

\bibitem[LGM11b]{Buzios}
Jean-Fran{\c{c}}ois Le~Gall and Gr\'egory Miermont.
\newblock {S}caling limits of random trees and planar maps.
\newblock ar{X}iv:1101.4856, 2011.

\bibitem[LGP08]{MR2438999}
Jean-Fran{\c{c}}ois Le~Gall and Fr{\'e}d{\'e}ric Paulin.
\newblock Scaling limits of bipartite planar maps are homeomorphic to the
  2-sphere.
\newblock {\em Geom. Funct. Anal.}, 18(3):893--918, 2008.

\bibitem[LL10]{MR2677157}
Gregory~F. Lawler and Vlada Limic.
\newblock {\em Random walk: a modern introduction}, volume 123 of {\em
  Cambridge Studies in Advanced Mathematics}.
\newblock Cambridge University Press, Cambridge, 2010.

\bibitem[LSW01a]{MR1879850}
Gregory~F. Lawler, Oded Schramm, and Wendelin Werner.
\newblock Values of {B}rownian intersection exponents. {I}. {H}alf-plane
  exponents.
\newblock {\em Acta Math.}, 187(2):237--273, 2001.

\bibitem[LSW01b]{MR1879851}
Gregory~F. Lawler, Oded Schramm, and Wendelin Werner.
\newblock Values of {B}rownian intersection exponents. {II}. {P}lane exponents.
\newblock {\em Acta Math.}, 187(2):275--308, 2001.

\bibitem[LW99]{MR1742883}
Gregory~F. Lawler and Wendelin Werner.
\newblock Intersection exponents for planar {B}rownian motion.
\newblock {\em Ann. Probab.}, 27(4):1601--1642, 1999.

\bibitem[Mie08]{MR2399286}
Gr{\'e}gory Miermont.
\newblock On the sphericity of scaling limits of random planar
  quadrangulations.
\newblock {\em Electron. Commun. Probab.}, 13:248--257, 2008.

\bibitem[Mie11]{Mu}
Gr\'egory Miermont.
\newblock {T}he {B}rownian map is the scaling limit of uniform random plane
  quadrangulations.
\newblock ar{X}iv:1104.1606, 2011.

\bibitem[Ons44]{MR0010315}
Lars Onsager.
\newblock Crystal statistics. {I}. {A} two-dimensional model with an
  order-disorder transition.
\newblock {\em Phys. Rev. (2)}, 65:117--149, 1944.

\bibitem[RV08]{arXiv:0807.1036}
R\'emi Rhodes and Vincent Vargas.
\newblock {KPZ} formula for log-infinitely divisible multifractal random
  measures.
\newblock arXiv:0807.1036, 2008.

\bibitem[RV10]{MR2642887}
Raoul Robert and Vincent Vargas.
\newblock Gaussian multiplicative chaos revisited.
\newblock {\em Ann. Probab.}, 38(2):605--631, 2010.

\bibitem[SD87]{MR889398}
H.~Saleur and B.~Duplantier.
\newblock Exact determination of the percolation hull exponent in two
  dimensions.
\newblock {\em Phys. Rev. Lett.}, 58(22):2325--2328, 1987.

\bibitem[She07]{MR2322706}
Scott Sheffield.
\newblock Gaussian free fields for mathematicians.
\newblock {\em Probab. Theory Related Fields}, 139(3-4):521--541, 2007.

\bibitem[She10]{ArXiv:1012.4797}
Scott Sheffield.
\newblock {C}onformal weldings of random surfaces: {SLE} and the quantum
  gravity zipper.
\newblock Ar{X}iv:1012.4797, 2010.

\bibitem[Smi10]{MR2680496}
Stanislav Smirnov.
\newblock Conformal invariance in random cluster models. {I}. {H}olomorphic
  fermions in the {I}sing model.
\newblock {\em Ann. of Math. (2)}, 172(2):1435--1467, 2010.

\end{thebibliography}
\bibliographystyle{alpha}

\end{document}